\documentclass{amsart}
\usepackage{latexsym,amsxtra,amscd,ifthen}
\usepackage{amsfonts}
\usepackage{verbatim}
\usepackage{amsmath}
\usepackage{amsthm}
\usepackage{amssymb}

\numberwithin{equation}{section}

\theoremstyle{plain}
\newtheorem{theorem}{Theorem}[section]
\newtheorem{lemma}[theorem]{Lemma}
\newtheorem{proposition}[theorem]{Proposition}
\newtheorem{corollary}[theorem]{Corollary}

\newcommand{\inth}{\textstyle \int}
\theoremstyle{definition}
\newtheorem{definition}[theorem]{Definition}
\newtheorem{example}[theorem]{Example}
\newtheorem{convention}[theorem]{Convention}
\newtheorem{remark}[theorem]{Remark}
\newtheorem{question}[theorem]{Question}

\makeatletter              
\let\c@equation\c@theorem  
                            
\makeatother

\DeclareMathOperator{\gldim}{gldim}

\DeclareMathOperator{\Ext}{Ext}

\DeclareMathOperator{\pdim}{projdim}

\DeclareMathOperator{\Spec}{Spec}

\DeclareMathOperator{\lcm}{lcm}

\DeclareMathOperator{\gr}{gr}

\DeclareMathOperator{\injdim}{injdim}

\DeclareMathOperator{\GKdim}{GKdim}
\DeclareMathOperator{\End}{End}

\DeclareMathOperator{\Hom}{Hom}
\DeclareMathOperator{\RHom}{RHom}

\DeclareMathOperator{\im}{im}

\DeclareMathOperator{\Mod}{{\sf Mod}}

\newcommand{\fm}{\mathfrak{m}}

\newcommand{\id}{\operatorname{id}}

\begin{document}

\title{Homological Integral of Hopf Algebras}

\author{D.-M. Lu, Q.-S. Wu and J.J. Zhang}

\address{Lu: Department of Mathematics, Zhejiang
University,
Hangzhou 310027, China}

\email{dmlu@zju.edu.cn}

\address{Wu: Institute of Mathematics, Fudan
University,
Shanghai, 200433, China}

\email{qswu@fudan.edu.cn}

\address{zhang: Department of Mathematics, Box 354350,
University of Washington, Seattle, Washington 98195,
USA}

\email{zhang@math.washington.edu}

\begin{abstract}
The left and right homological integrals are introduced 
for a large class of infinite dimensional Hopf algebras.
Using the homological integrals we prove a version of 
Maschke's theorem for infinite dimensional Hopf algebras. 
The generalization of Maschke's theorem and homological 
integrals are the keys to study noetherian regular Hopf 
algebras of Gelfand-Kirillov dimension one. 
\end{abstract}





\maketitle


\setcounter{section}{-1}
\section{Introduction}
\label{xxsec0}

Let $H$ be a Hopf algebra over a base field $k$. 
A result of Larson and Radford \cite{LR1,LR2} 
states that if $H$ is finite dimensional and if 
$\operatorname{char}k =0$, then the following two
conditions are
equivalent:

\smallskip

(GLD) $H$ has global dimension 0, namely, $H$ is
semisimple artinian.

(ANT) $S^2=\id_H$ where $S$ is the antipode of $H$.

\smallskip

\noindent
A result of Larson and Sweedler \cite{LS} states that
(GLD) is equivalent to:

\smallskip

(ITG) $\epsilon(\inth^r)\neq 0$ where $\inth^r$ is the
right
integral of $H$.

\smallskip
Larson and Sweedler's result (GLD) $\Leftrightarrow$ 
(ITG) is a generalization of Maschke's theorem 
for finite groups. 
These results are so elegant and useful that we can not help
attempting to extend them to the infinite dimensional case. The
extension of (ANT) $\Rightarrow$ (GLD) is quite successful. In
\cite[0.1]{WZ2} the authors proved the following: Suppose that 
$H$ is a finite module over its affine center and that 
$\operatorname{char}k =0$. If $S^2=\id_H$, 
then $H$ has finite global dimension. It is well-known that the 
converse of \cite[0.1]{WZ2} is false [Example \ref{xxex2.7}]. 
Until now the extension of (GLD) $\Leftrightarrow$ (ITG) is 
less successful, partly due to the fact that we don't have 
a good definition of the left and right integrals of an 
infinite dimensional Hopf algebra. In this paper we use 
homological properties to define the integrals for a large 
class of infinite dimensional Hopf algebras; and then 
generalize the result of Larson and Sweedler (GLD) 
$\Leftrightarrow$ (ITG).

To define the homological integral we need to assume that
$H$ is Artin-Schelter Gorenstein, which is satisfied by 
many noetherian Hopf algebras [Section \ref{xxsec1}]. 
The {\it left homological integral} $\inth^l$ of $H$ is the 
1-dimensional $H$-bimodule $\Ext^d_H({_Hk},{_HH})$ where $d$ 
is the injective dimension of $H$.  The right homological
integral $\inth^r$ is defined similarly [Definition 
\ref{xxdefn1.1}]. In the finite 
dimensional case the left and right homological integrals
agree with the usual left and right integrals of $H$ 
respectively. 

We say $H$ satisfies (Cond1) if the map
$$\Ext^d_H(\inth^r,\epsilon):\quad
\Ext^d_{H}(\inth^r,{_HH})\to \Ext^d_{H}(\inth^r,{_Hk})$$ 
is an isomorphism where $d$ is the injective dimension of 
$H$. We say $H$ satisfies (Cond2) if, for every simple 
left $H$-module $T\not\cong \inth^r$, $\Ext^d_H(T,k)=0$.
A $k$-algebra $A$ is called {\it regular} if $\gldim
A<\infty$ and called {\it affine} if it is finitely generated
as a $k$-algebra. Algebras satisfying a polynomial 
identity are called {\it PI algebras}. Our first theorem is 
a generalization of Larson and Sweedler's result (GLD) 
$\Leftrightarrow$ (ITG), which is a homological version of 
Maschke's theorem in the infinite dimensional setting.

\begin{theorem}
\label{xxthm0.1}
Suppose $H$ is a noetherian affine PI Hopf algebra.
Then $H$ is regular if and only if conditions
\textup{(Cond1)} and \textup{(Cond2)} hold.
\end{theorem}

Theorem \ref{xxthm0.1} follows from Theorems 
\ref{xxthm3.3} and \ref{xxthm3.4}. 
(Cond1) is not as pretty as (ITG); but when $H$ is 
finite dimensional, namely, when $d=0$,  (Cond1) is 
equivalent to (ITG) [Lemma \ref{xxlem3.5}].
(Cond2) is even more strange; but when $H$ is
finite dimensional, it is an easy consequence of
(Cond1) [Lemma \ref{xxlem3.5}]. Hence Theorem 
\ref{xxthm0.1} generalizes the result of Larson 
and Sweedler (GLD) $\Leftrightarrow$ (ITG). 
When $H$ is commutative, it is also easy to
see that (Cond2) is a consequence of (Cond1). 
Theorem \ref{xxthm0.1} would become much nicer 
if one can show that (Cond2) is a consequence 
of (Cond1) in general.

The definition of homological integrals uses only the
algebra structure of $H$, but not the coproduct 
of $H$. Any change of the coproduct of $H$ will not 
effect the homological integrals [Example \ref{xxex2.7}]. 
It is possible that there is another version of homological 
integrals that reflects the coalgebra structure of $H$.

In the finite dimensional case the left and right
integrals live in the $H$. In the infinite dimensional 
case we will define the residue module $\Omega$ of $H$ 
where both the left and the right homological integrals 
live. We say $H$ is {\it unimodular} if $\inth^l=\inth^r$ 
in $\Omega$. We define the {\it integral order} of $H$, 
denoted by $io(H)$, to be the minimal positive integer 
$n$ such that $(\inth^r)^{\otimes n}\cong k$ as 
$H$-bimodule [Definition \ref{xxdefn2.2}]. Then $H$ is 
unimodular if and only if $io(H)=1$. 

Integrals have been playing an important role in  
the studies of finite dimensional Hopf algebras. We 
expect that homological integrals will be useful in 
research into infinite dimensional Hopf algebras
of low Gelfand-Kirillov dimensions. The 
Gelfand-Kirillov dimension is denoted by 
GK-dimension from now on. In the second half 
of the paper we use homological integrals to 
investigate the structure of regular Hopf algebras 
of GK-dimension one. Using the right homological 
integral of $H$ we can define a quotient Hopf algebra
$H_{iq}$ of $H$, called the {\it integral quotient} 
of $H$ [Definition \ref{xxdefn4.2}], which plays 
an important role in the study of regular Hopf 
algebras of GK-dimension one.

\begin{theorem}[Theorem \ref{xxthm7.1}]
\label{xxthm0.2} Let $H$ be a noetherian affine Hopf
algebra of GK-dimension one and $H_{iq}$ be the 
integral quotient of $H$. Suppose $H$ is regular and prime. 
Then the following hold. 
\begin{enumerate}
\item
$io(H)(=\dim_k H_{iq})=PI.deg(H)$, where $PI.deg$ is the PI degree. 
\item
The coinvariant subalgebra $H^{co\; H_{iq}}$ is 
an affine commutative domain. 
\item
If $H$ is unimodular (or equivalently $H_{iq}$ is trivial), 
then $H$ is commutative. If further $k$ is algebraically closed, then
$H$ is isomorphic to either $k[x]$ or $k[x^{\pm 1}]$.
\end{enumerate}
\end{theorem}

Partial results and conjectural descriptions are given
when $H$ is not prime. Roughly speaking, every regular 
Hopf algebra $H$ of GK-dimension 1 should fit into a short 
exact sequence:
\begin{equation}
\label{E0.2.1}
0\to H_{dis}\to H\to H_{conn}\to 0
\tag{E0.2.1}
\end{equation}
where the connected component $H_{conn}$ is a regular
prime Hopf quotient algebra of $H$ and the discrete
component $H_{dis}$ {\it is conjecturally} a finite 
dimensional subalgebra of $H$ [Theorem \ref{xxthm6.5} 
and Remark \ref{xxrem6.6}]. The connected component
$H_{conn}$ should fit into a short exact sequence:
\begin{equation}
\label{E0.2.2}
0\to H_{cl}\to H_{conn} \to H_{iq}\to 0
\tag{E0.2.2}
\end{equation}
where the classical component $H_{cl}$ is the commutative
subalgebra $H_{conn}^{co\; H_{iq}}$ (see 
Theorem \ref{xxthm0.2}(b)) and the integral quotient 
$H_{iq}$ is the dual of a finite group algebra 
acting on $H_{conn}$ [Theorem \ref{xxthm7.1}
and Remark \ref{xxrem7.2}]. The statements can be found 
in Sections 6 and 7. 

Note that the conjectural descriptions of \eqref{E0.2.1} 
and \eqref{E0.2.2}  are verified for group algebras. 
According to \eqref{E0.2.1} and \eqref{E0.2.2} every 
group $G$ of linear growth is isomorphic to $G_{dis}
\rtimes({\mathbb Z}\rtimes G_{iq})$ [Proposition \ref{xxprop8.2}]. 
This description of groups of linear growth  
is well-known \cite{IS,St}. The ideas used here should 
be useful in further studies of homological properties 
and classification of infinite dimensional noetherian 
Hopf algebras and group algebras. For example we are 
wondering how the homological integrals can be used 
to study noetherian regular Hopf algebras of GK-dimension 
two. One of the testing questions is to extend 
Theorem \ref{xxthm0.2} to the GK-dimension two case
(see Examples \ref{xxex2.9} and \ref{xxex8.5}).

Using a slight generalization of Theorem 
\ref{xxthm0.1} (see Theorem \ref{xxthm3.3}(b)) and 
Theorem \ref{xxthm0.2}  we have the following easy 
corollary, which has a representation-theoretic flavor. 

\begin{corollary}
\label{xxcor0.3} Let $H$ be as in Theorem \ref{xxthm0.2}.
If $M$ is a simple left $H$-module, then $\dim_k M$
divides $io(H)$.
\end{corollary}

Corollary \ref{xxcor0.3} is proved at the end of 
Section 7. Note that Corollary \ref{xxcor0.3}
fails for the GK-dimension two case [Examples 
\ref{xxex2.9} and \ref{xxex8.5}].

Homological methods are effective for a large class of
infinite dimensional Hopf algebras, in particular, 
for noetherian affine PI Hopf algebras. We refer to 
\cite{Bro,BG,WZ1,WZ2} for some known results and for 
questions concerning the homological properties of 
these Hopf algebras.

\section{Definition of Integrals}
\label{xxsec1}

From now on let $k$ be a base field. There is no 
further restriction on $k$ unless otherwise stated.
We refer to Montgomery's book \cite{Mo} for the basic definitions 
about Hopf algebras. Let $H$ be a Hopf algebra over 
$k$. Usually $k$ denotes the trivial algebra or the
trivial Hopf algebra. For simplicity $k$ also 
denotes the trivial $H$-bimodule $H/\ker \epsilon$ 
where $\epsilon: H\to k$ is the counit of $H$. 
Usually we are working on left modules. Let 
$H^{\sf op}$ denote the opposite 
ring of $H$. A right $H$-module can be viewed as 
a left $H^{\sf op}$-module. An $H$-bimodule is 
sometimes identified with a left $H\otimes 
H^{\sf op}$-module where $\otimes$ denotes 
$\otimes_k$.

Recall that $H$ is {\it Artin-Schelter Gorenstein}
(or {\it AS-Gorenstein}) if

\smallskip
(AS1) $\injdim {_HH}=d<\infty$,

(AS2) $\dim_k \Ext^d_H({_Hk},{_HH})=1$,
$\Ext^i_H({_Hk},{_HH})=0$
for all $i\neq d$,

(AS3) the right $H$-module versions of the conditions
(AS1,AS2) hold.

\smallskip

We say $H$ is {\it Artin-Schelter regular} (or {\it
AS-regular}) if it is AS-Gorenstein and it has finite
global dimension.

It follows from the proof of \cite[Lemma 1.11]{BG} that
(AS2) implies

(AS2)${}'$ for each finite dimensional simple left $H$-module $M$,
$\dim \Ext^d_H(M,H)=\dim M$ and $\Ext^i_H(M,H)=0$ for all $i\neq
d$; and the same holds for right modules.

\begin{definition}
\label{xxdefn1.1}
Let $H$ be an AS-Gorenstein Hopf algebra of
injective dimension $d$. Any nonzero element in
$\Ext^d_H({_Hk},{_HH})$
is called a {\it left homological integral} of $H$. We
write
$\inth^l=\Ext^d_H({_Hk},{_HH})$. Any nonzero element
in
$\Ext^d_{H^{\sf op}}
({k_H},{H_H})$ is called a {\it right homological
integral} of $H$. We write
$\inth^r=\Ext^d_{H^{\sf op}}({k_H},{H_H})$. By abusing
the language we also call $\inth^l$ and $\inth^r$
the left and the right homological integrals of $H$
respectively.
\end{definition}

Homological integrals exist only for AS-Gorenstein Hopf algebras.
Hence free Hopf algebras (of at least two variables) and universal
enveloping algebras of infinite dimensional Lie algebras do not
have homological integrals. On the other hand we expect that
noetherian Hopf algebras are AS-Gorenstein, whence homological
integrals exist. It is well-known that finite dimensional Hopf
algebras are AS-Gorenstein of injective dimension 0. Affine
noetherian PI Hopf algebras are AS-Gorenstein by \cite[Theorem
0.1]{WZ1}. Many Hopf algebras associated to classical and quantum
groups are AS-Gorenstein and AS-regular \cite{Bro,BG}.

When $H$ is finite dimensional, then homological
integrals agree with the classical integrals
\cite[Definition 2.1.1]{Mo} in the
following way: the (classical) left integral is an
$H$-subbimodule of $H$; and it is identified with
the left homological integral $\Hom_{H}(k,H)$ via the
natural homomorphism
$$\Hom_H(\epsilon,H): \Hom_H(k,H)\to \Hom_H(H,H)\cong H.$$
The same holds for the right integral. 

Note that both $\inth^l$ and $\inth^r$ are 1-dimensional
$H$-bimodules. As a left $H$-module, $\inth^l\cong k$,
but as a right $H$-module, $\inth^l$ may not be isomorphic
to $k$. A similar comment applies to $\inth^r$.

\begin{definition}
\label{xxdefn1.2}
Let $H$ be a Hopf algebra with $\inth^l$ and
$\inth^r$.  We say $H$
is  {\it unimodular} if $\inth^l$ is isomorphic to $k$
as $H$-bimodules.
\end{definition}

The unimodular property means that
$$hx=xh=\epsilon(h)x$$
for all $h\in H$ and $x\in \inth^l$. When $H$ is finite
dimensional, this definition agrees with the classical definition
in \cite[p.~17]{Mo}.

\begin{lemma}
\label{xxlem1.3}
Suppose $H$ is noetherian. The following are
equivalent:
\begin{enumerate}
\item
$H$ is unimodular.
\item
$\inth^r\cong k$ as $H$-bimodules.
\item
$\inth^l\cong \inth^r$ as $H$-bimodules.
\end{enumerate}
\end{lemma}

\begin{proof} We only need to show that (a)
$\Leftrightarrow$ (b).
Since $H$ is noetherian and has finite injective
dimension, there is a convergent spectral sequence
\cite[(3.8.1)]{SZ}
$$\Ext^p_H(\Ext^q_{H^{\sf op}}({k_H},{H_H}),{_HH})
\Longrightarrow k_H.$$ 
The AS-Gorenstein condition gives rise to the 
isomorphism
\begin{equation}
\label{E1.3.1} \Ext^d_H(\Ext^d_{H^{\sf
op}}({k_H},{H_H}),{_HH})\cong k_H \tag{E1.3.1}
\end{equation}
where $d$ is the injective dimension of $H$.
If (b) holds, then
$\inth^r:=\Ext^d_{H^{\sf op}}({k_H},{H_H})\cong k$
as $H$-bimodule. Hence \eqref{E1.3.1} implies that
$\inth^l:=\Ext^d_H({_Hk},{_HH})\cong k$
as $H$-bimodule, whence $H$ is unimodular. By the
left-right symmetry we have the other implication.
\end{proof}

When $H$ is finite dimensional, $\inth^l$ and
$\inth^r$ live in the same
space $H$; and $H$ is unimodular if and only if
$\inth^l=\inth^r$.
This is also true for homological integrals 
after we define residue module of $H$.

An $H$-module $M$ is {\it locally finite} if every
finitely
generated submodule of $M$ is finite dimensional over
$k$.
An $H$-bimodule is {\it locally finite} if it is
locally finite on
both sides. Let $\Mod_{fd} H$ (respectively,
$\Mod_{fd} H^{\sf op}$)
denote the category of finite dimensional left
(respectively, right)
$H$-modules.

\begin{definition}
\label{xxdefn1.4}
Let $H$ be AS-Gorenstein of injective dimension $d$.
An $H$-bimodule
$\Omega$ is called a {\it residual module} of $H$ if
the following
conditions hold.
\begin{enumerate}
\item
$\Omega$ is locally finite.
\item
The functors $\Ext^d_H(-,{_HH})$ and
$\Hom_H(-,{_H\Omega})$ are naturally isomorphic 
when restricted to $\Mod_{fd}H \to
\Mod_{fd}H^{\sf op}$.
\item
The functors $\Ext^d_{H^{\sf op}}(-,{H_H})$ and
$\Hom_{H^{\sf op}} (-,{\Omega_H})$ are naturally 
isomorphic when restricted to
$\Mod_{fd} H^{\sf op}\to \Mod_{fd} H$.
\end{enumerate}
\end{definition}

Since $\Omega$ is locally finite, the above 
definition implies that $\Omega$ is unique up to 
a bimodule isomorphism. When $H$ is finite 
dimensional (namely, $d=0$), then $\Omega$ is 
isomorphic to $H$. The terminology ``residue
module'' indicates that $\Omega$ is related to
Yekutieli's residue complex \cite{YZ}

\begin{lemma}
\label{xxlem1.5}
Suppose $\Omega$ is the residual module of $H$.
\begin{enumerate}
\item
The left and right homological integrals
$\inth^l$ and $\inth^r$ can be identified
with $H$-subbimodule $\Hom_H({_Hk},{_H\Omega})$ and
$\Hom_{H^{\sf op}}(k_H,\Omega_H)$ respectively. With this
identification
$\inth^l$ and $\inth^r$ live in a common vector space
$\Omega$.
\item
$H$ is unimodular if and only if $\inth^l=\inth^r$.
\end{enumerate}
\end{lemma}

\begin{proof} (a) It follows from the natural
$H$-bimodule homomorphisms
$$\inth^l:=\Ext^d_H(k,H)\cong \Hom_H(k,\Omega)\to
\Hom_H(H,\Omega)=\Omega$$
that $\inth^l$ is a subbimodule of $\Omega$.
The same argument works for $\inth^r$.

(b) If $H$ is unimodular, then $\inth^l=kx\cong k$ as
$H$-bimodule,
for some $x\in \Omega$. Since $\inth^r\subset \Omega$
is the
1-dimensional
$H$-subbimodule isomorphic to $k$ as right $H$-module,
$\inth^r=kx=\inth^l$.
The converse is trivial.
\end{proof}

The existence of $\Omega$ is not all clear. We present
some partial results here.

Let $H$ be an affine noetherian PI Hopf algebra. By
\cite[Theorem 0.2(4)]{WZ1}, there is an exact
$H$-bimodule complex
\begin{equation}
\label{E1.5.1} 0\to H\to I^{-d}\to \cdots \to 
I^{-1}\to I^0\to 0
\tag{E1.5.1}
\end{equation}
such that \eqref{E1.5.1} is a minimal injective 
resolution of the left $H$-module $_HH$ and of 
the right $H$-module $H_H$ respectively. Further, 
when restricted the left or the right, $I^{-i}$ is 
pure of GK-dimension $i$. The complex \eqref{E1.5.1} 
is called {\it the residual complex} of $H$ by 
Yekutieli \cite{YZ}. In the present paper we 
are mainly interested in the last term $I^0$. As a left 
(respectively, right) $H$-module, $I^0$ is a union of 
injective hulls of finite dimensional left
(respectively, right) $H$-modules, which is locally
finite. In general, a residual complex exists for 
every noetherian affine PI AS-Gorenstein algebra 
\cite[Theorem 4.10]{YZ}.

\begin{lemma}
\label{xxlem1.6}
Suppose $H$ is a noetherian AS-Gorenstein Hopf algebra.
If $H$ has a residue complex \eqref{E1.5.1}, then $I^0$ 
is the residual module of $H$. In particular, the 
residual module $\Omega$ exists for every noetherian 
affine PI Hopf algebra.
\end{lemma}

\begin{proof}
Since the GK-dimension of $I^0$ is zero, it is locally
finite. Another property of the residual complex
is that $I^{-i}$ has no nonzero submodule of 
GK-dimension dimension less than $i$. If $i\neq 0$,
then $\Hom_H(M,I^{-i}) =0$ for all finite dimensional
left $H$-module $M$. Thus
$$\Hom_H(M,\eqref{E1.5.1})=\Hom(M,I^0[-d]),$$
which implies that
$$\Ext^d_H(M,{_HH})=\Hom(M,I^0).$$
Similarly,
$\Ext^d_{H^{\sf op}}(N,H_H)=\Hom_{H^{\sf op}}(N,I^0_H)$ 
for all finite dimensional right $H$-modules $N$. By 
definition $I^0$ is the residue module.
\end{proof}

Residual complexes may not exist for non-PI Hopf
algebras \cite[Remark 5.14]{YZ}. But the residue 
module could still exist in that case. We mention 
one result without proof: the residual module
exists for the enveloping algebra $U(\mathfrak g)$ of
a finite dimensional Lie algebra  $\mathfrak g$.

\begin{example}
\label{xxex1.7} Let $\mathfrak g$ be the simple Lie
algebra $sl_2$ generated by $e,f,h$ subject to relations
$$[e,f]=h,\quad [h,e]=2e,\quad [h,f]=-2f.$$
Let $H$ be the enveloping algebra $U(\mathfrak g)$.
Then $H$ is affine, noetherian and  AS-regular with
GK-dimension and global dimension 3.

It is well-known that $H$ has only one 1-dimensional
simple module, which is the trivial module. Hence $\inth^l=
\inth^r=k$ and $H$ is unimodular. It is not hard to verified
that $\Ext^i_H(k,k)=0$ for all $i\neq 0,3$ and $\Ext^i_H(k,k)
\cong k$ where $i=0,3$.

When $\operatorname{char} k>0$, $H$ is PI and $\Omega$
is isomorphic to $\lim (H/I)^*$ where $I$ runs over 
all co-finite dimensional ideals of $H$. Here $(-)^*$ 
denotes the $k$-linear vector space dual. Incidently,
$H^\circ$ is always defined to be $\lim (H/I)^*$. 
Here $H^\circ$ denotes the Hopf algebra dual in the
sense of \cite[Chapter 9]{Mo}.

When $\operatorname{char} k=0$, $H$ is not PI and $\Omega$
is isomorphic to $\oplus_{n\geq 0} M_n(k)$, which is also
equal to $\lim (H/I)^*$ where $I$ runs over all co-finite
dimensional ideals of $H$.

In general one can prove that $U(\mathfrak g)$ is 
unimodular if $\mathfrak g$ is a semisimple Lie algebra. 

By Example \ref{xxex3.2}, $U(\mathfrak g)$ may not 
be unimodular, if $\mathfrak g$ is not semisimple.
\end{example}

\section{Order of Integral}
\label{xxsec2}

Since homological integrals agree with usual 
integrals in finite dimensional case, we often 
use ``integral'' instead of ``homological integral'' 
from now on.

Let $\inth^l$ and $\inth^r$ be the left and the right 
integrals of $H$. For any $H$-bimodule $M$ let $S(M)$ 
denote the $H$-bimodule defined by the action
$$h'\cdot m\cdot h=S(h)mS(h')$$
for all $m\in M,h,h'\in H$. We can also define $S$ for 
one-sided $H$-modules in a similar way.

\begin{lemma}
\label{xxlem2.1} Let $H$ be a Hopf algebra with
integrals. Suppose the antipode $S$ 
of $H$ is bijective. Then $S(\inth^r)=\inth^l$ 
and $S(\inth^l)=\inth^r$.
\end{lemma}

\begin{proof} Since $S$ is bijective, the functor $S$
is invertible. For any $H$ bi-modules $M, N$, it is 
easy to check that
$$\Hom_H(S(M),S(N)) \cong S(\Hom_{H^{\sf op}}
(M,N))$$ 
as $H$ bi-modules. This isomorphism can be extended
to their derived functors. Obviously, 
$S(k)\cong k$  and the map $S:H\to H$ induces 
an isomorphism $H \cong S(H)$ as right and left 
and bi $H$-module. Hence we have
$$\Ext^d_H(k,H) \cong \Ext^d_{H}(S(k),S(H))) 
\cong S(\Ext^d_{H^{\sf op}}(k,H)).$$ 
Thus we proved $S(\inth^r)=\inth^l$. The proof of
$S(\inth^l)=\inth^r$ is the same.
\end{proof}

Given two left $H$-modules $M$ and $N$ we can 
define a left $H$-module structure on $M\otimes N$ 
via the coproduct $\Delta: H\to H\otimes H$. By 
the coassociativity of $\Delta$ the $n$th 
tensor product $M^{\otimes n}$ is a well-defined left 
$H$-module. Similarly we can define $\otimes$ for 
right $H$-modules. It is easy to check that 
$S(M\otimes N) \cong S(N) \otimes S(M)$ as right 
$H$-modules.

\begin{definition}
\label{xxdefn2.2} The integral order of $H$, 
denoted by $io(H)$, is the order of the right 
integral $\inth^r$, namely, the minimal 
positive integer $n$ (or $\infty$ if no such $n$) 
such that $(\inth^r)^{\otimes n}\cong k$ as 
left $H$-modules.
\end{definition}

The integral order has been used in the study 
of finite dimensional Hopf algebras either explicitly
or implicitly. 
By the above definition $io(H)$ seems dependent 
on the coproduct of $H$, because $\otimes$ is 
dependent on the coproduct of $H$.
By Lemma \ref{xxlem2.1} if the antipode $S$ is 
bijective, the $io(H)$ can also be computed by 
using the left integral.

Any $1$-dimensional left $H$-module 
$M$ can be identified as a quotient of an 
algebra homomorphism $\pi: H\to H/l.ann_H(M)$. 
Such an algebra homomorphism $\pi$ represents 
a group-like element in the dual Hopf algebra 
$H^{\circ}$. It is well-known that $\pi$ also 
defines an algebra automorphism $\sigma_{\pi}: 
H\to H$ given by
$$\sigma_{\pi}: h\mapsto \sum h_1\pi(h_2),$$
whose inverse is defined by
$$\sigma_{\pi}^{-1}: h\mapsto \sum h_1\pi(S(h_2)).$$
Let $\sigma^r$ be the automorphism of $H$ induced 
by the map $\Sigma^r: H\to H/l.ann(\inth^r)$. The 
following is clear.

\begin{lemma}
\label{xxlem2.3}
Suppose integrals exist for $H$. 
\begin{enumerate}
\item $io(H)$ is equal to the order of the 
group-like element $\Sigma^r \in H^{\circ}$. 
\item $io(H)$ is equal to the order of the
algebra automorphism $\sigma^r$ of $H$.
\end{enumerate}
\end{lemma}

Next we want to investigate $io(H)$ when $H$ 
is a noetherian PI ring.

Let $A$ be an algebra (not necessarily a Hopf 
algebra). We recall the definition of a 
clique in the prime spectrum $\Spec A$ and more 
details can be found in \cite[Chapter 11]{GW}. 
Let $P$ and $Q$ be two prime ideals. If there 
is a nonzero $A$-bimodule $M$ that is a 
subquotient of $(P\cap Q)/PQ$ and is torsionfree 
as left $A/P$-module and as right $A/Q$-module, 
then we say there is a link from $P$ to $Q$, 
written $P\leadsto Q$. The links make $\Spec A$ 
into a directed graph (i.e., a quiver) and the 
connected components of this graph are called
{\it cliques}.

We will only work with noetherian affine PI 
algebras $A$. Let $M$ and $N$ be two simple 
left $A$-modules. We say $M$ and $N$ are in 
the same clique if $l.ann_A(M)$ and 
$l.ann_A(N)$ are in the same clique. It is 
easy to check that $l.ann_A(M)\leadsto l.ann_A(N)$ 
if and only if $\Ext^1_A(N,M)\neq 0$ (this
is a consequence of \cite[Theorem 11.2]{GW}, 
see also \cite[p.~324]{BW}). In this case we 
sometimes write $M\leadsto N$. The following 
lemma is more or less known.

\begin{lemma}
\label{xxlem2.4} Let $A$ be a noetherian affine 
PI algebra. Let $M$ and $N$ be two 
simple left $A$-modules. If $\Ext^n_A(M,N)\neq
0$ for some $n$, then $M$ and $N$ are in the 
same clique.
\end{lemma}

\begin{proof} We prove a slightly more general 
statement and the assertion follows from the 
general statement.

{\it Claim:} Let $X$ be any locally finite 
left $A$-module. Let
$$0\to X\to I^0\to I^1\to \cdots \to I^n\to \cdots$$
be the minimal injective resolution of $X$. 
Then every simple subquotient of $I^n$ is in 
the same clique as some simple subquotient of $X$.

By induction we may assume that $n=0$. Let $S$
be a simple subquotient of $I^0$. Pick an 
element $x\in I^0$ such that the submodule 
$Y$ generated by $x$ has a quotient module
isomorphic to $S$. If $Y\cap X$ is not in 
the kernel of $Y\to S$, then $S$ is a 
quotient of $Y\cap X$. 
We are done. Otherwise, $Y/Y\cap X \to S$ is
surjective. In this case by the induction
on $\dim_k Y/Y\cap X$ we can assume that
$\dim_k Y/Y\cap X=1$ (note that we can also 
change $X$ when using this induction). Since 
$Y$ is an essential extension of $Y\cap X$, 
we have $\Ext^1_A(S,Y\cap X)\neq 0$. This 
implies that there is a simple subquotient 
$S'$ of $Y\cap X$ such that $\Ext^1_A(S,S')
\neq 0$. Hence $S$ and $S'$ are linked and 
we have proved our {\it claim}.

Now let $X$ be the simple module $N$. If 
$\Ext^n_A(M,N)\neq 0$, then $M$ is a 
subquotient of $I^n$ where $I^n$ is the
$n$th term in the minimal injective 
resolution of $N$. By the {\it claim} 
$M$ and $N$ are in the same clique.
\end{proof}

\begin{proposition}
\label{xxprop2.5} Let $H$ be an affine 
noetherian PI Hopf algebra. Suppose that 
$\Ext^i_{H}(\inth^r,k)\neq 0$ for some
$i$. If the clique containing $k$ is finite, 
then $io(H)$ is finite.
\end{proposition}

\begin{proof} Since $\inth^r$ is 1-dimensional,
by \cite[Lemma 1.3]{WZ2} 
\begin{equation}
\label{E2.5.1}
(S(\inth^r))^*\otimes \inth^r\cong k\cong 
\inth^r\otimes (S(\inth^r))^*
\tag{E2.5.1}
\end{equation} 
as left $H$-modules. Note that our 
definition of $^*$ is slightly different 
from the definition given in \cite[p.~602]{WZ2} 
and that's the reason we need to add $S$. 
The equation \eqref{E2.5.1} implies that
the functor $\inth^r\otimes -$ is an 
auto-equivalence of the the category of left 
$H$-modules. Hence 
$\Ext^i_{H}(M,N)\cong \Ext^i_{H}(\inth^r\otimes M,
\inth^r\otimes N)$ for all left $H$-modules
$M,N$ and for all $i$. In particular,
$$\Ext^i_{H}((\inth^r)^{\otimes (p+1)},
(\inth^r)^{\otimes p})
\cong \Ext^i_{H}(\inth^r,k)\neq 0$$
for all $p$. By Lemma \ref{xxlem2.4}, all 
$(\inth^r)^{\otimes p}$'s are in the clique 
containing $k$. Since that clique is finite
by hypotheses, 
there are $n>m$ such that 
$(\inth^r)^{\otimes n}\cong (\inth^r)^{\otimes m}$.
Hence $(\inth^r)^{\otimes (n-m)}\cong k$ and
$io(H)\leq n-m$.
\end{proof}

\begin{lemma}
\label{xxlem2.6} Let $H$ be an AS-Gorenstein 
Hopf algebra and let $x$ be a normal 
non-zero-divisor of $H$ such that $(x)$ is a 
Hopf ideal of $H$. Suppose that $\tau$ is the
algebra automorphism of $H$ such that 
$xh=\tau(h)x$ for all $h\in H$. 
\begin{enumerate}
\item
$H':=H/(x)$ is an AS-Gorenstein Hopf algebra.
\item
$\inth^l_H\cong (\inth^l_{H'})^{\tau^{-1}}$ as
right $H$-modules.
\item
If $x$ is central, then 
$\inth^l_H \cong \inth^l_{H'}$ and $io(H)=io(H')$.
\end{enumerate}
\end{lemma}

\begin{proof} Let $M$ be a left $H'$-module. 
By change of rings \cite[Theorem 11.66 or 
Corollary 11.68]{Ro}, we have 
$$\Ext^p_H(M,H)\cong \Ext^p_H(H/(x)\otimes_{H/(x)} M,H)
\cong \Ext^{p-1}_{H/(x)}(M,\Ext^1_H(H/(x),H)).$$
An easy computation shows that $\Ext^1_H(H/(x),H)
\cong {^\tau}(H/(x))\cong (H/(x))^{\tau^{-1}}$. Hence
$$\Ext^p_H(M,H)\cong\Ext^{p-1}_{H/(x)}(M,(H/(x))^{\tau^{-1}})
\cong \Ext^{p-1}_{H/(x)}(M,H/(x))^{\tau^{-1}}.$$
This shows that $\injdim H'\leq \injdim H-1$.
Let $M=k$. Then we see that $\injdim H'=\injdim H-1$
and $\inth^l_H=(\inth^l_{H'})^{\tau^{-1}}$ as right
$H$-modules. The AS-Gorenstein property follows 
from this  isomorphism. Finally when $x$ is 
central, then $\tau$ is the identity map of 
$H$ and $\inth^l_H=\inth^l_{H'}$. In this 
case $io(H)=io(H')$.
\end{proof}

We will have a finiteness result about $io(H)$ 
in Lemma \ref{xxlem5.3}(g). To conclude this 
section we give two examples, the first of 
which comes from Taft's construction \cite{Ta} 
(see also \cite[Example 1.5.6]{Mo}).

\begin{example}
\label{xxex2.7}
Let $n$, $m$ and $t$ be integers and $\xi$ be an
$n$-th primitive root of $1$. Let $H$ be the 
$k$-algebra generated by $x$ and $g$ subject to 
the relations 
$$g^n=1, \quad\text{and}\quad xg=\xi^m gx.$$
So $H$ is commutative if and only if $\xi^m=1$.
The coalgebra structure of $H$ is defined by
$$\Delta(g)=g\otimes g, \; \epsilon(g)=1,
\quad\text{and}\quad \Delta(x)=x\otimes 
1+g^t\otimes x, \; \epsilon(x)=0.$$
So $H$ is cocommutative if and only if $g^t=1$. 
The antipode $S$ of $H$ is defined by
$$S(g)=g^{-1}\quad\text{and}\quad
S(x)=-g^{-t}x=-\xi^{mt}xg^{-t}.$$
The order of the antipode $S$ is 
$2\; order(\xi^{mt})$. Hence $H$ 
is involutive if and only if $\xi^{mt}=1$. 
One way to understand the algebra structure 
of $H$ is to view it as a smash product 
$k[x]\# kG$ where $G$ is the group of 
$\langle g\rangle$. The $kG$ action on $k[x]$ 
is determined by $$g\circ x=\xi^m x.$$
When $\operatorname{char} k$ does not divide 
$n$, by \cite[Corollary 2.4]{LL}, the global 
dimension of $H$ is bounded by the global 
dimension of $k[x]$. Hence $\gldim H=1$. 
Since $x^n$ is a central element, $H$ is 
finite over its center. Therefore $H$ is 
AS-regular. Change $t$ to another $t'$ will 
only change the coproduct of $H$, but not 
the algebra structure of $H$. So this 
change will not effect the (homological) 
integrals of $H$. To compute the integrals 
we use Lemma \ref{xxlem2.6}. Note that 
$x$ is a normal element with $xh=\tau(h)x$ for 
all $h\in H$ and where $\tau: x\to x, g\to 
\xi^m g$. Since $H'=H/(x)$ is isomorphic 
$kG$, it commutative and finite
dimensional. Hence $\inth^l_{H'}=k=H/(x,g-1)$.
By Lemma \ref{xxlem2.6}, $\inth^l_H
=H/(x,g-1)^{\tau^{-1}}\cong H/(x,g-\xi^{-m})$ as
right $H$-modules. By Lemma \ref{xxlem2.1},
$\inth^r\cong H/(x,g-\xi^{m})$ as left $H$-modules.

The integral order $io(H)$ is equal to the
$order(\xi^m)$. When $\gcd(m,n)=1$, then
$H$ is a prime ring of PI degree $n$. In this
case the PI degree of $H$ is equal to $io(H)$.
In this example, any change of coproduct 
of $H$ does not effect $io(H)$ though the 
definition of $io(H)$ uses the coproduct of $H$.

If $\xi^{mt}=-1$, then the ideal generated by $x^2$
is a Hopf ideal. The quotient Hopf algebra is not
semisimple because the order of $S$ is 4. So a 
quotient Hopf algebra of a regular algebra may 
not be regular. On the other hand, certain 
quotient Hopf algebras of a regular algebras are
regular if they have the same GK-dimension, 
see Lemma \ref{xxlem5.5}(a).
\end{example}

\begin{lemma}
\label{xxlem2.8} Let $H$ and $K$ be two noetherian 
AS-Gorenstein Hopf algebra such that $H\otimes K$ 
is noetherian with finite injective dimension. Then 
\begin{enumerate}
\item $H\otimes K$ is AS-Gorenstein and 
$\injdim H\otimes K= \injdim H+\injdim K$. 
\item
$\inth^r_{H\otimes K}=\inth^r_H\otimes \inth^r_K$
and 
$io(H\otimes K)=\lcm \{io(H),io(K)\}$.
\end{enumerate}
\end{lemma}

\begin{proof} Since $k_{H\otimes K}=k_{H}\otimes k_{K}$, 
we have 
$$\Ext^i_{H\otimes K}(k,H\otimes K)
\cong \bigoplus_{j=0}^i \Ext^j_H(k,H)\otimes 
\Ext^{i-j}_K(k,K).$$ The assertions follow from
the AS-Gorenstein property of $H$ and $K$.
\end{proof}

\begin{example}
\label{xxex2.9}
Let $H$ be the Hopf algebra in Example \ref{xxex2.7}
with $m=t=1$ and $n>1$. Then $H$ is prime and 
$io(H)=PI.deg(H)=n$. It is easy 
to check that $H\otimes H$ is prime, noetherian, affine
PI of injective dimension 2 and that $PI.deg(H\otimes H)=n^2$. 
By Lemma \ref{xxlem2.8}, $io(H\otimes H)=n$. Since
$PI.deg(H\otimes H)=n^2$, there is a simple 
$H\otimes H$-module of dimension $n^2$. Therefore
Corollary \ref{xxcor0.3} fails for Hopf algebras 
of GK-dimension 2. 

In general, $io(K)\neq PI.deg(K)$ for a nice
Hopf algebra $K$ of GK-dimension $>1$. Theorem 
\ref{xxthm0.2} fails badly for higher GK-dimension.
\end{example}

\section{Extension of Larson-Sweedler}
\label{xxsec3}

As mentioned in the introduction, Larson and 
Sweedler proved the following version of Maschke's 
theorem: a finite dimensional Hopf 
algebra is semisimple artinian (i.e., has
global dimension 0) if and only if 
$\epsilon(\inth^l)\neq 0$, and if and only if 
$\epsilon(\inth^r)\neq 0$. The terms 
$\epsilon(\inth^l)$ and $\epsilon(\inth^r)$ 
make sense in the infinite dimensional case
in the following way.
The counit $\epsilon: H\to k$ induces an $H$-bimodule
homomorphism, which is also denoted by $\epsilon$,
\begin{equation}
\label{E3.0.1} \epsilon: \inth^l
=\Ext^d_H({_Hk},{_HH})\to
\Ext^d_H({_Hk},{_Hk}).\tag{E3.0.1}
\end{equation}

\begin{lemma}
\label{xxlem3.1}
If $\epsilon(\inth^l)\neq 0$, then $H$ is unimodular.
\end{lemma}

\begin{proof}
Note that $\Ext^d_H({_Hk},{_Hk})$ is isomorphic 
to a direct sum of $k$ as $H$-bimodules. If
$\epsilon(\inth^l)\neq 0$, then $\epsilon$ in 
\eqref{E3.0.1} embeds $\inth^l$ into  
$\Ext^d_H({_Hk},{_Hk})$. Hence $\inth^l$ is 
isomorphic to $k$.
\end{proof}

To generalize Larson-Sweedler's result it is 
natural to ask if the condition 
$\epsilon(\inth^l)\neq 0$ in Lemma \ref{xxlem3.1}
is equivalent to $H$ having finite global 
dimension. The answer is ``No'' as the next example 
shows (see also Example \ref{xxex2.7}).

\begin{example}
\label{xxex3.2}
Let $L$ be the 2-dimensional solvable Lie algebra
generated by $x$ and $y$ subject to the relation
$[x,y]=x$. Let $H$ be the enveloping algebra 
$U(L)$. Then $H$ is a noetherian affine domain 
of global dimension 2. It is involutory, i.e.,
$S^2=1$. If $\operatorname{char} k\neq 0$, 
then it is a PI algebra.

The trivial module $k$ is isomorphic to $H/(x,y)$. 
By a computation in \cite[Example 3.2]{ASZ} or
using Lemma \ref{xxlem2.6} we have
$$\inth^l:=\Ext^2_H({_Hk},{_HH})\cong
H/(x,y+1)\not\cong k$$
as right $H$-modules. For every $p$, 
$(\inth^l)^{\otimes p}\cong H/(x,y+p)$ as right
$H$-modules. Hence $io(H)=\operatorname{char} k$
if $\operatorname{char} k\neq 0$ or $io(H)=\infty$
otherwise. In particular, $H$ is not unimodular.  
By Lemma \ref{xxlem3.1}, $\epsilon(\inth^l)=0$.
\end{example}

Conditions (Cond1) and (Cond2) are the correct 
replacement of (ITG) as the next result shows. Recall 
from the introduction that $H$ satisfies (Cond1) 
if the map
$$\Ext^d_H(\inth^r,\epsilon):\quad
\Ext^d_{H}(\inth^r,{_HH})\to
\Ext^d_{H}(\inth^r,k)$$
is an isomorphism where $d$ is the injective
dimension of $H$ and that $H$ satisfies (Cond2) if
for every simple left $H$-module $T\not\cong
\inth^r$, $\Ext^d_H(T,k)=0$.

\begin{theorem}
\label{xxthm3.3}
Suppose $H$ is noetherian and AS-regular of 
global dimension $d$.
\begin{enumerate}
\item 
\textup{(Cond1)} and \textup{(Cond2)} hold. 
\item 
Let $W$ be a finite dimensional simple left 
$H$-module and let $W'$ denote the simple left 
$H$-module $\Ext^d_H(W,H)^*$. Then, for 
every surjective map $H\to W'$ of left $H$-modules, 
the induced map $\Ext^d_H(W,H)\to \Ext^d_H(W,W')$ 
is nonzero and surjective. Further, for every 
simple left $H$-module $T\not\cong W$,
$\Ext^d_H(T,W')=0$.
\end{enumerate}
\end{theorem}

\begin{proof} If $W=\inth^r$, then $W'=
\Ext^d_H(\inth^r,H)^* \cong (k_H)^*\cong {}_Hk$ 
(see Lemma \ref{xxlem1.3} for the argument of 
$\Ext^d_H(\inth^r,H)\cong k_H$). Also 
$\Ext^d_H(\inth^r,H)$ is 1-dimensional, the 
assertion (b) implies (a). So we only prove (b).

By \cite[Proposition 7.1]{ASZ}, for every finite
dimensional $H$-module $W$ and every noetherian 
left $H$-module $N$ there is a natural isomorphism
\begin{equation}
\label{E3.3.2}
\Ext^{d-i}_H(W,N)\cong \Ext^i_H(N,W')^*
\tag{E3.3.2}
\end{equation}
where $W'=\Ext^d_H(W,H)^*$ and where $*$ is the
$k$-linear dual. For $i=0$,
the functor $\Ext^d_H(W,-)$ is equivalent to
$\Hom_H(-,W')^*$.
This equivalence translates the assertions in (b) 
into the following assertions about $\Hom_H(-,W')$:

(i) the induced map $\Hom_H(H,W')^*\to \Hom_H(W',W')^*
(\cong k)$ is surjective; and

(ii) $\Hom_H(W,T)=0$ for all simple $T\not\cong W$.

\noindent
Both (i) and (ii) are obviously true now.
\end{proof}

Theorem \ref{xxthm3.3} gives one implication in 
Theorem \ref{xxthm0.1}. In the rest of this 
section we prove the other implication in Theorem 
\ref{xxthm0.1}. We introduce a condition slightly 
weaker than those in Theorem \ref{xxthm3.3}. We 
say $H$ satisfies (Cond3) if, for every simple 
left $H$-module $T$, 
$$\Ext^d_H(T,\epsilon):
\Ext^d_H(T,{_HH})\to \Ext^d_H(T,{_Hk})$$ 
is surjective. It is obvious that (Cond3) is 
a consequence of (Cond1)+(Cond2)
(when taking $W=\inth^r$ and $W''=k$). The 
following is a converse of Theorem \ref{xxthm3.3}.

\begin{theorem}
\label{xxthm3.4}
Let $H$ be noetherian and AS-Gorenstein.
Suppose that injective hulls of finite dimensional
left $H$-modules are locally finite. If
\textup{{(Cond3)}} holds, then $H$ is regular.
\end{theorem}

Combining Theorem \ref{xxthm3.3} and \ref{xxthm3.4},
we see that $H$ is regular if and only if \textup{{(Cond3)}}
holds.

\begin{proof}[Proof of Theorem \ref{xxthm3.4}] 
Let $d$ be the injective dimension of
$H$. First we prove the following two statements:
\begin{enumerate}
\item
$\Ext^d_H(-,{_Hk})$ is a right exact functor on
$\Mod_{fd} H$.
\item
$\Ext^d_H(M,{_HH})\to \Ext^d_H(M,{_Hk})$ is surjective
for all $M\in \Mod_{fd}H$.
\end{enumerate}

\noindent
Consider an exact sequence in $\Mod_{fd}H$,
\begin{equation}
\label{E3.4.1}
0\to L\to M\to N\to 0 \tag{E3.4.1}
\end{equation}
with nonzero $L$ and $N$. Applying $\Ext^d_H(-,H)$ 
and $\Ext^d_H(-,k)$ to \eqref{E3.4.1} we obtain
a row exact commutative diagram
\begin{equation}
\label{E3.4.2}
\CD
0 @>>> \Ext^d_H(N,H) @>>> \Ext^d_H(M,H) @>>>
\Ext^d_H(L,H) @>>> 0\\
@VVV     @VfVV @VgVV @VhVV @VVV \\
K @>>> \Ext^d_H(N,k) @>p>> \Ext^d_H(M,k) @>e>>
\Ext^d_H(L,k) @>>> C\\
\endCD
\tag{E3.4.2}
\end{equation}
for some modules $K$ and $C$. The top row is 
exact since $H$ is AS-Gorenstein. The bottom
row is exact for $K=\ker (p)$ and $C=
\operatorname{coker}(e)$.

We first prove (b) by induction on $\dim_k M$. By 
(Cond3) the assertion (b) holds when $M$ is simple 
(in particular when $\dim_k M=1$). If $\dim_k M>1$ 
and $M$ is not simple then we can find nonzero $L$ 
and $N$ that fit into \eqref{E3.4.1}. By induction
hypothesis the vertical maps $f$ and $h$ are
surjective. By a diagram chasing or by 
\cite[Lemma 3.32 (five lemma)]{Ro}, $g$ is 
surjective. Then (b) follows from induction.

Now assume (b). To prove (a) we need to show 
that the map $e$ is surjective. This is clear 
since $g$ and $h$ are surjective by (b).

Secondly we show that (a) implies that $_Hk$ has
finite injective dimension. Let
\begin{equation}
\label{E3.4.3}
0\to {_Hk}\to J^0\to J^1\to \cdots \to J^d\to
J^{d+1}\to  \cdots
\tag{E3.4.3}
\end{equation}
be a minimal injective resolution of $_Hk$.

Suppose that $k$ has infinite injective dimension.
Then $J^{d+1}\neq 0$. Let $T=\ker(J^d\to J^{d+1})$.
Since \eqref{E3.4.3} is minimal, $T$ is an essential proper
submodule of $J^d$. Let 
$\Phi$ be the set of all co-finite dimensional 
ideals of $H$. For every $I\in \Phi$, let
$$T_I=\Hom_H(H/I,T)\quad \text{and}\quad
J_I=\Hom_H(H/I,J^d).$$ Since $H$ is noetherian, 
$\Ext^d_H(M,k)$ is finite dimensional over $k$ 
for every finitely generated left $H$-module $M$. 
When $M$ is finite dimensional and simple, $J^d$
contains only finitely many copies of $M$ for each such 
$M$. For a fixed $I\in \Phi$, the socle of $J_I$ 
is finite dimensional since there are finitely 
many simple modules over $H/I$. This implies that
$J_I$ is finite dimensional. Since $J^d$ is 
locally finite, we have
$$T=\lim_{I\in \Phi} T_I\quad \text{and}\quad
J^d=\lim_{I\in \Phi} J_I.$$
Then the fact $T\neq J^d$ implies that $T_I\neq J_I$ 
for some $I$. Now we consider the exact sequence 
of finite dimensional left $H$-modules
$$0\to T_I\to J_I\to N\to 0$$
where $N=J_I/T_I$ is nonzero. Because any
map $f:T_I\to T$ factoring through $J^{d-1}\to T$
can not be injective, the inclusion map
$\alpha: T_I\to T$ induces a nonzero element
$$\overline{\alpha}\in
\Ext^d_H(T_I,k)=\Hom_H(T_I,T)/\im(\Hom_H(T_I,J^{d-1})).$$
We claim that there is no element
$$\overline{\beta}\in
\Ext^d_H(J_I,k)=\Hom_H(J_I,T)/\im(\Hom_H(J_I,J^{d-1}))$$
such that
\begin{equation}
\label{E3.4.4}
\overline{\alpha}=\overline{\beta} \circ i
\tag{E3.4.4}
\end{equation}
where $i: T_I\to J_I$ is the inclusion. The equation
\eqref{E3.4.4} means that there is a map $\beta: 
\Hom_H(J_I,T)$ such that $\alpha=\beta \circ i+
\partial^{d-1}\circ \phi$ where $\phi$ is some 
element in $\Hom_H(T_I,J^{d-1})$ and $\partial^{d-1}:
J^{d-1}\to J^d$ is the differential map of the
complex \eqref{E3.4.3}. To show 
the claim we rewrite the equation as 
$\alpha-\partial^{d-1}\circ \phi=\beta\circ i$.
For any simple submodule $S$ of $T_I$, the image 
$\phi(S)$ must be in the socle of $J^{d-1}$. Hence
$\partial^{d-1}\circ \phi(S)=0$ since \eqref{E3.4.3} 
is minimal. This implies that the map $\alpha':=
\alpha-\partial^{d-1}\circ \phi$ is injective on
the socle of $T_I$, whence injective on $T_I$. Since
the socle of $T_I$ is equal to the socle of $J_I$ 
by minimality of \eqref{E3.4.3}, the equation 
$\alpha'=\beta\circ i$ implies that the map 
$\beta$ is injective on the socle of $J_I$, 
whence injective on $J_I$. But this is impossible 
because $\dim J_I>\dim T_I$. Therefore $\beta$ 
does not exist.

Non-existence of $\beta$ shows that the sequence
$$\to \Ext^d_H(N,k)\to \Ext^d_H(J_I,k)\to
\Ext^d_H(T_I,k)\to 0$$
is not right exact. This yields a contradiction with
(a). This contradiction shows that $_Hk$ has finite
injective dimension.

Since $H$ is noetherian and has finite injective 
dimension, $H$ itself serves as a dualizing 
complex over $H$. By \cite[Lemma 2.1]{WZ3}, the 
right $H$-module $\inth^l:=\Ext^d_H(k,H)\cong 
\RHom_H(k,H)[d]$  (where $[d]$ is the $d$th 
complex shift) has finite right projective 
dimension. Since $\inth^l$ is 1-dimensional, by
\cite[Proposition 1.4]{WZ2}, $\gldim H$ is finite.
\end{proof}

Theorem \ref{xxthm0.1} follows from Theorems 
\ref{xxthm3.3} and \ref{xxthm3.4} because 
injective hulls of finite dimensional modules 
over noetherian affine PI rings are locally finite.

Suppose $H$ is commutative. Then for any two 
non-isomorphic simple modules $S_i$, we have 
$\Ext^i_H(S_1,S_2)=0$. This implies that 
$\Ext^i_H(S,k)=0$ for all $i$ all $S\not
\cong k$. Since $H$ is commutative, $\inth^r=\inth^l=k$. 
Hence (Cond2) is automatic. A similar statement holds 
for finite dimensional Hopf algebras.

\begin{lemma}
\label{xxlem3.5}
Suppose $H$ is finite dimensional. Then the following
are
equivalent:
\begin{enumerate}
\item
\textup{(ITG)}.
\item
\textup{(Cond1)}.
\item
\textup{(Cond1)} plus \textup{(Cond2)}.
\end{enumerate}
\end{lemma}

\begin{proof} It's easy to see that (ITG) implies that
$H$ is unimodular. The same applies to (Cond1). So we
might as well assume $H$ is unimodular. Under this hypothesis,
(Cond2) becomes trivial; and (ITG) is just (Cond1).
\end{proof}

Based on the limited evidences in the finite 
dimensional case and in the commutative case we 
ask the following question.

\begin{question}
\label{xxque3.6}
Let $H$ be a noetherian affine PI Hopf algebra of
injective dimension $d$. If  \textup{(Cond1)} 
holds, is then $H$ regular?
\end{question}

\section{Integral Quotient of $H$}
\label{xxsec4}

In this section we discuss some quotient
(i.e., factor) Hopf 
algebras of $H$ which will be used later. The
following is clear.

\begin{lemma}
\label{xxlem4.1} Let $H$ be a Hopf algebra and
let $\fm$ be an ideal of $H$. Then there is
a unique maximal Hopf ideal $J$ contained in $\fm$.
\end{lemma}

\begin{definition}
\label{xxdefn4.2} Let $\inth^r$ be the right 
integral of $H$. Let $\fm=l.ann_H(\inth^r)$
and let $J_{iq}$ be the maximal Hopf ideal contained in
$\fm$. The quotient Hopf algebra $H/J_{iq}$ is called 
the {\it integral quotient} of $H$, denoted by 
$H_{iq}$.
\end{definition}

Another Hopf quotient related to $H_{iq}$ is the 
abelianization of $H$. The following lemma is clear.

\begin{lemma}
\label{xxlem4.3}
Let $H$ be a Hopf algebra. Let $I$ be the ideal generated
by $xy-yx$ for all $x,y\in H$. Then $I$ is a Hopf ideal
and $H\to H/I$ is a Hopf algebra homomorphism.
\end{lemma}

We call $H/I$ in the above lemma the {\it abelianization}
of $H$, denoted by $H_{ab}$. Let $\pi$ be the 
canonical Hopf algebra homomorphism $H\to H_{ab}$.
If $f: H\to H'$ is a Hopf algebra homomorphism 
with $H'$ being commutative, then $f$ factors through 
$\pi: H\to H_{ab}$. As in Example \ref{xxex7.3} 
the map $\pi: H\to H_{ab}$ is not necessarily conormal 
in the sense of \cite[Definition 3.4.5]{Mo}. Since 
$\fm:=l.ann_H(\inth^r)$ has co-dimension 1, $I\subset
\fm$. This implies that $I\subset J$. Therefore 
$H_{iq}$ is a Hopf quotient of $H_{ab}$. As a consequence
$H_{iq}$ is a commutative Hopf algebra.

\begin{lemma} 
\label{xxlem4.4}
Let $H$ be a Hopf algebra with integrals.
\begin{enumerate}
\item
$io(H)$ is finite if and only if $\dim_k H_{iq}$ is
finite.
\item
If $io(H)=n<\infty$, then $H_{iq}\cong (k {\mathbb Z}_n)^\circ$.
\item
$J_{iq}=\bigcap_{i=0}^{n-1}\ker ((\Sigma^r)^i: H\to k)$.
\end{enumerate}
\end{lemma}

\begin{proof} Suppose first $H_{iq}$ is finite dimensional.
Let $\Sigma^r: H\to H/\fm$ where $\fm=l.ann_H(\inth^r)$. 
Then $\Sigma^r$ is a group-like element in $H_{iq}^\circ
\subset H^\circ$. Hence $\bigoplus_{i\in {\mathbb Z}} k(\Sigma^r)^i$
is a Hopf subalgebra of $H_{iq}^\circ$.  By 
Lemma \ref{xxlem2.3}(a) $io(H)$ is bounded by 
$\dim_k H_{iq}^\circ=\dim H_{iq}$. 

Now assume that $io(H)=n<\infty$. Then the finite
group algebra 
$\bigoplus_{i=0}^{n-1} k(\Sigma^r)^i$ is a Hopf subalgebra 
of $H_{iq}^\circ$. Hence $I:=\bigcap_{i=0}^{n-1}
\ker ((\Sigma^r)^i: H\to k)$ is a Hopf ideal of $H$.
Since $\ker (\Sigma^r: H\to k)=l.ann_H(\inth^r)=\fm$ and
since $\fm$ contains the Hopf ideal $J_{iq}$, by
a computation using coproduct, each 
$\ker ((\Sigma^r)^i: H\to k)$ contains $J_{iq}$.
Hence we have $I\supset J_{iq}$. By the maximality of 
$J_{iq}$, we have $I=J_{iq}$. Thus $\dim_k H_{iq}=n$. 
This forces that $H_{iq}^\circ =\bigoplus_{i=0}^{n-1} 
k(\Sigma^r)^i$. The assertions follow.
\end{proof}

Let $G$ be a group. Define $G_{ab}$ to be $G/[G,G]$.
The following lemma holds because $H_{iq}$ is a 
Hopf quotient of $H_{ab}$. 

\begin{lemma}
\label{xxlem4.5} Suppose integrals exist for $H$.
If $H_{ab}$ is finite dimensional, then $io(H)$ divides 
$\dim H_{ab}$. If $H$ is a group algebra $kG$, then $io(H)$ 
divides $|G_{ab}|$.
\end{lemma}

The next example shows that $H_{iq}\neq H_{ab}$ in
general even when both are finite dimensional.

\begin{example}
\label{xxex4.6} 
Let ${\mathbb D}$ denote the group $\langle g,x
| g^2=1, gxg=x^{-1}\rangle$. It is easy to
see that ${\mathbb D}$ contains ${\mathbb Z}$
(generated by $x$) as a normal subgroup and 
${\mathbb D}/{\mathbb Z}\cong {\mathbb Z}_2$.
Also ${\mathbb D}$ is isomorphic to the free
product ${\mathbb Z}_2*{\mathbb Z}_2$. The
Hopf algebra $k {\mathbb D}$ (when 
$\operatorname{char}k\neq 2$) is prime and regular 
of GK-dimension 1. The abelianization $H_{ab}$ is 
$H/(x^2-1)$, and isomorphic to 
$k({\mathbb Z}_2\oplus {\mathbb Z}_2)$. 

To compute the integral, we note that
$x-1$ is a normal element of $H$ since
$$(x-1)g=(-gx^{-1})(x-1).$$
Using this fact and Lemma \ref{xxlem2.6} one sees that
$\inth^r\cong H/(x-1,g+1)$ as left $H$-modules. Hence 
$io(H)=2$ and $H_{iq}=H/(x-1)\cong k\; {\mathbb Z}_2$.
Therefore $H_{ab}\neq H_{iq}$.
\end{example}

We will also use the following easy lemma.

\begin{lemma}
\label{xxlem4.7}
Let $I$ be any Hopf ideal of $H$. Let $J=\bigcap_n I^n$.
Then $J$ is a Hopf ideal of $H$.
\end{lemma}

\section{Noetherian PI Hopf algebras}
\label{xxsec5}

The first step beyond the finite dimensional case is to look 
into affine noetherian Hopf algebras of GK-dimension 1. 
By \cite{SSW}, such an algebra is PI. In this
section we collect some known results about noetherian
affine PI Hopf algebras. In later sections we will
concentrate on Hopf algebras of GK-dimension one.

\begin{lemma}
\label{xxlem5.1}
Let $H$ be a noetherian affine PI  Hopf algebra of
GK-dimension $d$.
\begin{enumerate}
\item
$H$ is AS-Gorenstein of injective dimension d. As a
consequence, the left and right integrals 
of $H$ exist.
\item
$H$ has a quasi-Frobenius ring of fractions, denoted
by $Q$.
\item
The residue module of $H$ exists.
\item
If $d=1$, the injective module $Q/H$ is the
residue module of $H$.
\end{enumerate}
\end{lemma}

\begin{proof} (a) This is \cite[Theorem 0.1]{WZ1}.

(b) This is \cite[Theorem 0.2(2)]{WZ1}.

(c,d) This follows from \cite[Theorem 0.2(4)]{WZ1} and
Lemma 1.6.
\end{proof}

The following lemma is \cite[Theorem 0.1 and Corollary 3.7]{WZ2}.

\begin{lemma}
\label{xxlem5.2} Let $H$ be a noetherian affine PI Hopf
algebra. Suppose that the base field $k$ is of
characteristic zero and that $H$ is involutive.
If either $H$ is finite over its affine center or
$\GKdim H=1$, then $H$ is regular.
\end{lemma}

The following lemma is a consequence of
\cite[Theorems 5.4, 5.6 and Remark 5.7]{SZ}. The
definition of a hom-hom PI ring given in
\cite[p.~1013]{SZ} is equivalent to our
definition of AS-regular. Note that a Krull domain
of dimension 1 is a Dedekind domain
(of global dimension 1). We refer to \cite{SZ, WZ1,WZ2}
for the definitions of Auslander regular and 
Cohen-Macaulay. 

\begin{lemma}
\label{xxlem5.3}
Let $H$ be a noetherian regular affine PI Hopf algebra.
\begin{enumerate}
\item
$H$ is AS regular, Auslander regular and Cohen-Macaulay
and $\GKdim H=\gldim H$.
\item
$H$ is a direct sum of prime rings of the same GK-dimension
and the center of $H$ is a direct sum of Krull domains of the
same GK-dimension.
\item
$H$ is finite over its center and each prime direct summand
is equal to its trace ring.
\item
Every clique of $H$ is finite and localizable.
\item
If $\GKdim H=1$, then the center of $H$ is a direct
sum of Dedekind domains.
\item
Let $P_0$ be the minimal prime ideal of $H$ contained in
$\ker \epsilon$ and let $H_0=H/P_0$. Then the number of maximal
ideals in the clique containing $k$ is bounded by the
PI degree of $H_0$.
\item
$io(H)\leq PI.deg(H_0)<\infty$.
\end{enumerate}
\end{lemma}

\begin{proof} By \cite[Theorem 1.14]{BG} $H$ is AS-regular.

(a-e) follow from \cite[Theorems 5.4 and 5.6]{SZ}.
We need to show (f) and (g).

(f) By (b) we have $H=\oplus_{i=0}^s H_i$ where each
$H_i$ is a prime component of $H$ and $\oplus_{i>0}H_i$
is contained in $\ker \epsilon$. Then $P_0=\oplus_{i>0}H_i$.
By (c) $H_0$ equals to its trace ring. The assertion follows
from \cite[Theorem 8]{Bra}.

(g) The integral is only dependent on the
homological property of $H_0$. By forgetting the coalgebra
structure of $H$ we might assume that $H$ is prime.
By \textup{(Cond1)}, $\Ext^d_H(\inth^r,k)\neq 0$ where
$d=\GKdim H$. By the proof of Proposition 2.5, $io(H)$
is bounded by the number of primes in the clique containing
$k$. The assertion follows from (f).
\end{proof}

\begin{definition}
\label{xxdefn5.4}
Let $H$ be a noetherian Hopf algebra with 
finite GK-dimension. Suppose $\pi: H\to H_0$ is a Hopf 
algebra quotient map. We say $H_0$ is the {\it connected
component} of $H$ if
\begin{enumerate}
\item 
$\GKdim H_0=\GKdim H$;
\item 
The following universal property holds: for every 
quotient Hopf algebra homomorphism $f: H\to H'$ 
of the same GK-dimension, there is a unique Hopf 
algebra homomorphism $g: H'\to H_0$ such that 
$gf=\pi$.
\end{enumerate}
The connected component of $H$ is denoted by 
$H_{conn}$. We say $H$ is {\it connected} if 
$H$ is a connected component of itself.
\end{definition}

A connected component of $H$ is unique if it exists. 
So it is safe to call $H_{conn}$  {\it the} 
connected component of $H$. If $H$ is prime,
then it is connected. The converse is not clear to us, 
though we will show that this is true for regular 
Hopf algebras of GK-dimension one 
[Theorem \ref{xxthm6.5}]. We hope to see that this 
is true for noetherian affine PI Hopf algebras.

\begin{lemma}
\label{xxlem5.5}
Let $H$ be a noetherian affine PI regular Hopf algebra
and let $H'$ be a quotient Hopf algebra of $H$.
Suppose that $\GKdim H=\GKdim H'=d$.
\begin{enumerate}
\item
$H'$ is regular. Further $\inth^l_H=\inth^l_{H'}$
as right $H$-modules and $io(H')=io(H)$.
\item
$H'$ is a direct summand of $H$.
\item
$H_{conn}$ exists. 
\end{enumerate}
\end{lemma}

\begin{proof} By \cite[Theorem 0.3]{WZ1}, $H'$ is
projective over $H$ on both sides. Then we have, 
for all $i>d$ and all left $H'$-module $M$,
$$\Ext^i_{H'}({_{H'}k},{_{H'}M})=
\Ext^i_{H'}(H'\otimes_H k,{_{H'}M})\cong
\Ext^i_{H}(k,\Hom_{H'}(H',M))=0.$$
Hence $k$ has finite projective dimension
over $H'$. By \cite[Corollary 1.4(c)]{BG},
$H'$ is regular. By Lemma \ref{xxlem5.3}(b),
$H'$ is a direct sum of prime rings of the
same GK-dimension. Thus $H'$ is a factor ring
of $H$ modulo some prime components. Hence
$H\cong H'\oplus A$ for some algebra $A$.
Since the trivial $H$-module is also an
$H'$-module we have $\Ext^d_H(k,H)\cong
\Ext^d_{H'}(k,H')$. This shows that
$\inth^l_H\cong \inth^l_{H'}$ as right
$H$-modules. Similarly, $\inth^r_H=\inth^r_{H'}$.
This implies that $(\inth^l_H)^{\otimes p}\cong
(\inth^l_{H'})^{\otimes p}$ for all $p$
as right $H$-module. The assertion of
$io(H)=io(H')$ follows. 

Finally we need to prove (c), namely, to 
show that the universal property in Definition 
\ref{xxdefn5.4}. By Lemma \ref{xxlem5.3}(b), 
$H$ is a direct
sum of prime rings of the same GK-dimension.
So we can write it as $H=\oplus H_i$ as 
a decomposition of algebras. By the above proof,
we see that every quotient Hopf algebra $H'$
is the factor ring $H/P$ where $P=\oplus_{j}H_j$
where $\{j\}$ is a subset of $\{i\}$. If $P$ and $Q$
are such Hopf ideals, then $P+Q$ is also a Hopf
ideal. Since $P+Q\neq H$ because both $P$ and $Q$ are
contained in $\ker \epsilon$. Therefore $P+Q=
\oplus_j H_j$ is a direct summand of $H$. Let $M$ 
be the union of all such Hopf ideals. Then $H_0=H/M$ 
has the desired universal property and hence 
$H_0=H_{conn}$.
\end{proof}

\section{Regular Hopf Algebras of GK-dimension 1}
\label{xxsec6}

The statement in Lemma \ref{xxlem5.3}(b) asserting that
$H$ is a direct sum of primes is analogous to the
classical fact of algebraic groups: Suppose
$\operatorname{char} k=0$. Let $G$ be an algebraic
group and let $G_0$ be the {\it connected component}
of the identity of $G$. Then $G$ is a finite disjoint
union of $gG_0$ for some $g\in G$ and there is a
short exact sequence of algebraic groups
$$0\to G_0\to G\to G/G_0\to 0$$
where $G/G_0$ is a finite group, which is called the
{\it discrete part} of $G$. Dually, for
a noetherian regular commutative Hopf algebra $H$,
we have a short exact sequence of Hopf algebras
\begin{equation}
\label{E6.0.1}
0\to H_{dis}\to H\to H_{conn}\to 0
\tag{E6.0.1}
\end{equation}
where $H_{dis}$ is the maximal finite dimensional
normal Hopf subalgebra of $H$ and where $H_{conn}$
is the quotient Hopf algebra $H/(H_{dis})^{+}H$.
The Hopf algebra $H_{conn}$  is a domain.

In this section we are aiming for a similar
statement for noetherian regular noncommutative
Hopf algebras of GK-dimension 1. Due to the
noncommutativity of $H$ there are some extra
structures of $H$ related to the integrals
of $H$.  For simplicity throughout this section 

{\it let $H$ be a
noetherian affine PI regular Hopf $k$-algebra of 
GK-dimension 1.} 

\noindent
In some parts
we also assume that $k$ is algebraically closed as will be
stated, but this assumption is not needed in the main
result. Recall that $A^*$ means the $k$-linear dual of
$A$. 

Recall that the {\it integral quotient} of $H$,  
denoted by $H_{iq}$, is isomorphic to $H/J_{iq}$. 
When $is(H)=n<\infty$, $J_{iq}=\bigcap_{p=0}^{n-1}
\ker ((\Sigma^r)^p: H\to k)$ [Lemma \ref{xxlem4.4}]. 
Or equivalently, $J_{iq}=
\bigcap_p l.ann_H((\inth^r)^{\otimes p})$.
By Lemma \ref{xxlem4.4}, $H_{iq}^\circ\cong k
{\mathbb Z}_n$. Of course if $k$ is 
algebraically closed and $\operatorname{char} k
\nmid n$, then $H_{iq}\cong 
k{\mathbb Z}_n$. The definitions of normal Hopf 
subalgebras, normal Hopf ideals and conormal 
homomorphisms are given in \cite[pp.~33-36]{Mo}.

\begin{proposition}
\label{xxprop6.1} Let $J_{iq}=\bigcap_p l.ann_H
((\inth^r)^{\otimes p})$ and $H_{iq}=H/J_{iq}$. 
Then the Hopf algebra homomorphism $H_{ab}\to 
H_{iq}$ is conormal. 
\end{proposition}

We need a few lemmas to prove this proposition. 
First we need a reduction to the case when $k$ 
is algebraically closed.

\begin{lemma}
\label{xxlem6.2}
Let $F$ be a field extension of the base field $k$.
Let $H_F=H\otimes F$.
\begin{enumerate}
\item
$H_F$ is a noetherian
affine PI regular Hopf algebra of GK-dimension 1.
\item
$\inth^r_{H_F}=\inth^r\otimes F$ and
$\inth^l_{H_F}=\inth^l\otimes F$.
\item
$io(H)=io(H_F)$.
\item
$J_F:=
\bigcap_p l.ann_{H_F}((\inth^r_{H_F})^{\otimes p})
=\bigcap_p l.ann_H((\inth^r)^{\otimes p})\otimes F$.
As a consequence, $J_F$ is a Hopf ideal of $H_F$ (if
and only if $J_{iq}$ is a Hopf ideal of $H$).
\item
Proposition \ref{xxprop6.1} holds for $H$ if and
only if it holds for $H_F$.
\end{enumerate}
\end{lemma}

\begin{proof} (a) Field extension clearly
preserves the following properties:
``noetherian'', ``affine PI'' and ``GK-dimension 1''.
To prove that $H_F$ is regular, we only need to
show that the trivial $H_F$-module $F$ has
finite projective dimension. This follows from
the fact that
$$\pdim_{H_F} F=\pdim_{H\otimes F} (k\otimes F)=
\pdim_{H} k=1$$
when $H$ is noetherian.

(b,c) Since $H$ is noetherian, by K\"unneth formula,
$$\inth^l_{H_F}=\Ext^1_{H_F}(F,H_F)=
\Ext^1_{H\otimes F}(k\otimes F,H\otimes F)\cong
\Ext^1_{H}(k,H)\otimes F=\inth^l\otimes F.$$
Similarly, $\inth^r_{H_F}=\inth^r\otimes F$ and
$io(H_F)=io(H)$.

(d) This is also clear from (b).

(e) This is true because of (b) and the fact
the algebra $H/J_{iq}$ is isomorphic to a 
finite direct sum of $k$. 
\end{proof}

By the above lemma we only need to prove Proposition
\ref{xxprop6.1} for $H_F$ where $F$ is the
algebraic closure of $k$. In other words we may
assume $k$ to be algebraically closed. The following
lemma is obvious because $H_{iq}=k$.

\begin{lemma}
\label{xxlem6.3} If $H$ is unimodular, then
Proposition \ref{xxprop6.1} holds trivially.
\end{lemma}

Next we deal with the case when $io(H)\neq 1$.
Let $H^{\circ}$ be the dual Hopf algebra of $H$.
For any ideal $P$ of $H$ of codimension one, let
$\pi_P: H\to H/P\cong k$ denote the corresponding
group-like element in $H^{\circ}$. Let $M_P$ be
the left $H$-module $H/P$. We identify
the ideal $P$ with the module $M_P$ when talk about
the cliques. This is convenient for us as we did in
Section 2. Note that if $P$ and $Q$ are two
ideals of $H$ of codimension one, then there is 
another ideal $R$ of $H$ of codimension one
such that $M_P\otimes M_Q=M_R$ where
the module structure on the tensor 
is defined via the coproduct of $H$.
In this case it also implies $\pi_P\pi_Q=\pi_R$ in
$H^\circ$. Let $l.ann_H(\inth^r)=P_0$. Then
$\inth^r=M_{P_0}$ and, for every $t$,  
$(\inth^r)^{\otimes t}=(M_{P_0})^{\otimes t}$ is
corresponding to $(\Sigma^r)^t$ in $H^\circ$.

\begin{lemma}
\label{xxlem6.4} Let $Q$ be an ideal of $H$ of 
codimension one.
\begin{enumerate}
\item
$M_Q\otimes \inth^r\cong \inth^r\otimes M_Q$.
\item
The clique containing $M_Q$ is $\{M_Q\otimes
(\inth^r)^{\otimes t}\; |\; t\geq 0\}$.
\item
The group $\{(\Sigma^r)^{t}\}_{t\geq 1}$ is a
finite central subgroup of the group 
generated by $\{\pi_P\}\subset H^\circ$ where 
$P$ runs over all possible ideals of $H$ 
of codimension one. 
\end{enumerate}
\end{lemma}

\begin{proof} First we prove that the clique
containing
$\ker \epsilon$ is the set of primes associated to
$\{(\inth^r)^{\otimes t}\; |\; t\geq 0\}$. By
\cite[p.~324]{BW},
for any two maximal ideals $P$ and $Q$, 
$Q\leadsto P$ if and
only if $\Ext^1_H(H/P,H/Q)\neq 0$. By Theorem
\ref{xxthm3.3} and
(Cond1,2), $\Ext^1_H(\inth^r,k)\neq 0$ and
$\Ext^1_H(T,k)=0$
for all other simple $H$-modules $T$. Hence 
the only left $H$-module linked to $\inth^r$ is $k$. 
Since $M_P\otimes -$ and $-\otimes M_P$ are 
equivalences of the category of left $H$-modules,
the only module linked to $M_P\otimes \inth^r$ 
is $M_P$. For the same reason, the only module
linked to $\inth^r\otimes M_P$ is $M_P$. Therefore 
$M_P\otimes \inth^r=\inth^r\otimes M_P$. This 
proves (a). Let $P=P_0$. We obtain 
that $\{(\inth^r)^{\otimes t}\;|\; t\geq 0\}$ is the 
clique containing $k$. Again using the functor
$M_P\otimes -$, we obtain (b). By definition
$\{(\inth^r)^{\otimes t}\;|\; t\geq 0 \}$ is a
group, that is finite by Lemma \ref{xxlem5.3}(g).

(c) follows from (a). 
\end{proof}

\begin{proof}[Proof of Proposition \ref{xxprop6.1}]
If $H$ is unimodular, see Lemma \ref{xxlem6.3}.

Now assume $io(H)\neq 1$. By Lemma \ref{xxlem6.2},
we may assume that $k$ is algebraically closed.
We will show the following.
\begin{enumerate}
\item[(i)]
The abelianization $H_{ab}$ is finite dimensional
and semisimple.
\item[(ii)]
The group algebra $\oplus_p k (\Sigma^r)^{p}$
is a normal Hopf subalgebra of $(H_{ab})^\circ$.
\end{enumerate}
If $H_{ab}$ has GK-dimension 1, then,
by Lemma \ref{xxlem5.5}(a), $io(H)=io(H_{ab})=1$.
This contradicts our hypothesis of $io(H)\neq 1$.
Hence $H_{ab}$ is finite dimensional.
Since $k$ is algebraically closed of characteristic
zero, $H_{ab}$ is a dual of a finite group algebra.
This shows (i).

The dual Hopf algebra $(H_{ab})^\circ$ is the group
algebra $\oplus k\pi_{Q_i}$. By Lemma \ref{xxlem6.4}(c)
$\oplus_p k (\Sigma^r)^{p}$ is a normal Hopf 
subalgebra of $(H_{ab})^\circ$. This is (ii).

Note that (ii) is equivalent to that the map 
$H_{ab}\to (\oplus_p k (\Sigma^r)^{p})^\circ$ is conormal.
Proposition \ref{xxprop6.1} follows from the
fact that $H_{iq}=(\oplus_p k (\Sigma^r)^{p})^\circ$
[Lemma \ref{xxlem4.4}].
\end{proof}

Here is the main result of this section.

\begin{theorem}
\label{xxthm6.5} Let $H=\oplus_{i=0}^d H_{i}$
be a decomposition of $H$ into prime components
as in Lemma \ref{xxlem5.3}(b).
Let $P_0$ be a minimal prime ideal contained in
$\ker \epsilon$.
\begin{enumerate}
\item
$P_0$ is a Hopf ideal and $H/P_0$ is the connected
component
of $H$ in the sense of Definition \ref{xxdefn5.4}, denoted by
$H_{conn}$.
\item
The coinvariant subalgebra $H^{co H_{conn}}$, denoted
by $H_{dis}$, contains all finite dimensional normal
Hopf subalgebras of $H$.
\end{enumerate}
\end{theorem}

\begin{remark}
\label{xxrem6.6} 
Both $H_{conn}$ and $H_{dis}$ are canonical objects
associated to $H$. As is suggested by the commutative 
case and Theorem \ref{xxthm6.5}, it is natural to 
ask if $H_{dis}$ is finite dimensional. It is also 
unclear to us if  $H_{dis}$ is a normal Hopf 
subalgebra of $H$. If these questions have 
positive answers we conjecture that there is a 
short exact sequence of Hopf algebras in the 
sense of \cite[Definition 1.5]{Sc}
\begin{equation}
0\to H_{dis}\to H\to H_{conn}\to 0
\tag{E0.2.1}
\end{equation}
where $H_{conn}=H/(H_{dis})^{+}H$.
The conjectural description \eqref{E0.2.1} is an
analog of \eqref{E6.0.1} in the commutative case.
This will be verified when $H$ is a group 
algebra [Proposition \ref{xxprop8.2}].
Even if \eqref{E0.2.1} is false we can still ask if 
there is a Hopf algebra structure
on the algebra $H_{dis}$ such that $H$ is isomorphic
to a cross product $H_{dis}\#_{\sigma} H_{conn}$ 
for some $2$-cocycle $\sigma$.
\end{remark}

The point of Theorem
\ref{xxthm6.5}(a) is that $H_{conn}$ is a prime
algebra which is not clear from Lemma \ref{xxlem5.5}.
We need a few lemmas before proving Theorem
\ref{xxthm6.5}.

\begin{lemma}
\label{xxlem6.7} If there is a Hopf algebra quotient map
$H\to H'$ such that $H'$ is prime and $\GKdim H'=1$,
then $H'$ is isomorphic to $H_{conn}$.
\end{lemma}

\begin{proof} This follows from the existence of
$H_{conn}$ in Lemma \ref{xxlem5.5} and the universal property
of $H_{conn}$ stated in Definition \ref{xxdefn5.4}.
\end{proof}

\begin{lemma}
\label{xxlem6.8} Assume the notation as in
Theorem \ref{xxthm6.5}.
\begin{enumerate}
\item
If $K$ is a finite dimensional
normal Hopf subalgebra of $H$, then $H/K^{+}H$
is infinite dimensional.
\item
If $H$ is prime, then there is no non-trivial
finite dimensional normal Hopf subalgebra of
$H$.
\end{enumerate}
\end{lemma}

\begin{proof} 
(a) [Brown] Let $K$ be a finite dimensional
normal Hopf subalgebra of $H$. Then $H$ is a free
left $K$-module \cite[Theorem 2.1(2)]{Sc}, with 
a basis $\{h_i\;|\; i\in I\}$. Since $K$ is finite 
dimensional and $H$ is infinite dimensional, $I$
is infinite. Note that $H/K^{+}H$ is isomorphic to
$\oplus_{i\in I}kh_i$. Thus $H/K^{+}H$ is
infinite dimensional.

(b) Let $K$ be a finite dimensional normal Hopf
subalgebra of $H$. If $K\neq k$, then $\bar{H}
:= H/K^{+}H$ is finite dimensional because $H$ 
is prime of GK-dimension one, a contradiction to (a).
\end{proof}

\begin{lemma}
\label{xxlem6.9} Let $J_{iq}$ be 
as in Proposition \ref{xxprop6.1}. Then
$\bigcap_n J_{iq}^n=P_0$.
\end{lemma}

\begin{proof} First of all $l.ann(k)$ contains
the minimal prime $P_0$. Since $H$ is a direct
sum of primes and the ideals in the clique
containing $k$ are related by non-vanishing of 
$\Ext^1_H$, the ideals in the clique also contain 
$P_0$. Thus the assertion is equivalent to
$\bigcap_n I^n=0$ where $I$ is the
image of $J_{iq}$ in $H/P_0$. Replacing $H$ by $H/P_0$
(and forgetting the coalgebra structure of $H$ for
a moment), it suffices to show that $\bigcap_n I^n=0$.
By Lemma 5.3(d) $H$ is localizable at $I$ and the
localization, denoted by $Q$, is semilocal.
Since $H\to Q$ is injective, we only need to show that
$\bigcap I'^n=0$ where $I'=IQ$ is the Jacobson radical
of $Q$. This is true for noetherian semilocal PI
rings \cite[Corollary 6.4.15]{MR} and so we have proved 
the assertion.
\end{proof}

\begin{proof}[Proof of theorem \ref{xxthm6.5}]
Since $J_{iq}$ is a Hopf ideal
of $H$, by Lemma \ref{xxlem4.7}, $\bigcap_n J_{iq}^n$ 
is a Hopf ideal. By Lemma \ref{xxlem6.9}, $P_0$
is a Hopf ideal. Hence $H\to H/P_0$ is a Hopf
algebra quotient with $H/P_0$ being prime and
GK-dimension 1. By Lemma \ref{xxlem6.7}, $H/P_0$
is the connected component of $H$. Thus (a)
is proved.

(b) Suppose $K$ is a finite dimensional
normal Hopf subalgebra of $H$. Let $\bar{H}=H/K^{+}H$.
By \cite[Proposition 3.4.3]{Mo}, $H^{co \bar{H}}=K$.
By (a) and Lemma \ref{xxlem6.8}(a), there is a canonical Hopf
algebra quotient map from $\bar{H}\to H_{conn}$.
Hence $K\subset H_{dis}$.
\end{proof}

\begin{corollary}
\label{xxcor6.10}
Let $H$ be prime and $F$ be a
field extension of $k$. Then $H_F:= H\otimes F$
is prime.
\end{corollary}

\begin{proof} By Lemma \ref{xxlem6.9}, if $H$ is prime
then $\bigcap_n J_{iq}^n=0$. By Lemma \ref{xxlem6.2}(d),
$\bigcap_n J_{F}^n=0$ where $J_F=
\bigcap_p l.ann_{H_F} ((\inth^r_{H_F})^{\otimes p})$. 
By Lemma \ref{xxlem6.9}, $H_F$ is prime.
\end{proof}

\section{Regular Hopf Algebras of GK-dimension 1: Prime case}
\label{xxsec7}

Throughout Section 7 

{\it let $H$ be a noetherian affine PI regular 
Hopf algebra of GK-dimension 1.} 

\noindent
Here is the main result, which is a restatement of Theorem 
\ref{xxthm0.2}.

\begin{theorem}
\label{xxthm7.1} Let $H$ be prime.
\begin{enumerate}
\item If $H$ is unimodular, then $H$ is commutative.
If further $k$ is algebraically closed, then $H$ is
isomorphic to either $k[x]$ or $k[x^{\pm 1}]$.
\item If $H$ is not unimodular with $io(H)=n$, then
the following holds.
\begin{enumerate}
\item[(i)]
$n=PI.deg(H)$.
\item[(ii)]
The subalgebra of coinvariants $H_{cl}:=
H^{co H_{iq}}$ is a commutative domain.
\item[(iii)]
$H_{cl}$ is affine and $H$ is a finite module over $H_{cl}$
on both sides.
\end{enumerate}
\end{enumerate}
\end{theorem}

In a previous version of this paper there was an extra 
hypothesis ``$\operatorname{char} k\nmid n$'' 
in Theorem \ref{xxthm7.1}(b)(iii). We thank 
Brown for providing us a proof of the current version 
of the statement.

\begin{remark}
\label{xxrem7.2} 
\begin{enumerate}
\item 
The extra hypothesis of $k=\bar{k}$ 
in the second assertion of Theorem \ref{xxthm7.1}(a) 
is needed. When $k$ is not algebraically closed, there
are other interesting Hopf algebras 
[Example \ref{xxex8.3}]. 
\item
An immediate question after Theorem \ref{xxthm7.1} is
whether $H_{cl}$ is a Hopf subalgebra of $H$. This 
is not true as the next example shows. But we don't
know if $H_{cl}$ is always equipped with a Hopf
structure which may not be compatible with
the Hopf structure of $H$. If this is the case,
we can ask if there is a ``twisted'' short exact 
sequence
$$0\to H_{cl}\to H\to H_{iq}\to 0$$
of Hopf algebras and if $H$ is isomorphic to
a cross product $H_{cl}\#_{\sigma} H_{iq}$.
\item
Theorem \ref{xxthm7.1} fails for Hopf algebras
of GK-dimension 2 (see Examples \ref{xxex2.9} and \ref{xxex8.5}). 
\end{enumerate}
\end{remark}

\begin{example}[Continuation of Example \ref{xxex2.7}]
\label{xxex7.3} 
Let $H$ be the Hopf algebra defined in Example 
\ref{xxex2.7} for $n>1,m=1$. Since $m=1$, $H$ is
prime. When $\operatorname{char} k\nmid n$,
then $H$ is regular. It is easy to see that 
$H_{iq}=H/(x)\cong 
k\langle g\rangle \cong k{\mathbb Z}_n$ and 
$H_{cl}$ is the subalgebra generated by $x$. 
Since $\Delta(x)=x\otimes 1+g\otimes x$, $H_{cl}$
is not a subcoalgebra of $H$. The ideal
$\ker(H\to H_{iq})$ is not normal in the sense of
\cite[Definition 3.4.5]{Mo}, or the homomorphism
$H\to H_{iq}$ is not conormal. Let $H_1$ be the Hopf
algebra $k[x]$ with $\Delta(x)=x\otimes 1+1\otimes x$.
Then $H_1\cong H_{cl}$ as algebras. It is easy to
see that $H_1$ is an $H_{iq}$-Hopf module algebra
and $H$ is a smash product $H_1\# k{\mathbb Z}_n$
[Example \ref{xxex2.7}]. This does suggest a ``twisted''
short exact sequence
$$0\to H_1\to H\to k {\mathbb Z}_n\to 0.$$
Also in this example, one has $H_{iq}=H_{ab}$, while
the group algebra in Example \ref{xxex4.6} does not
have this property. 

Note that $\fm:=\ker \epsilon$ is idempotent.
Hence $H_{\fm}^\circ$ is isomorphic to the trivial
Hopf algebra $k$. On the other hand, $J_{iq}$ satisfies
the condition $\bigcap_n J^n_{iq}=0$ [Lemma \ref{xxlem6.9}].
Hence $H^\circ_{J_{iq}}$ is dense in $H^*$ (see the
proof of \cite[Proposition 9.2.10]{Mo}).
\end{example}

We need a few lemmas before proving Theorem
\ref{xxthm7.1}. A module is called {\it uniserial} 
if its submodules form a chain under inclusion. We 
will consider uniserial modules with a simple 
submodule. Let $T$ be a finite dimensional simple 
left $H$-module. Let $T'$ denote the left 
$H$-module $\Ext^1_H(T,H)^*$. Since both 
$\Ext^1_H(-,H)$ and $(-)^*$ are equivalent
functors on finite dimensional $H$-modules, $T'$ is
necessarily finite dimensional simple and $\End_H(T)
=\End_H(T')$. Since $H$ is a Hopf algebra it follows
from (AS2)$'$ in Section 1 that $\dim T'=\dim T$. Let $E(T)$ be the
injective hull of $T$. Since $H$ is PI and $T$ is finite
dimensional, $E(T)$ is locally finite.

\begin{lemma}
\label{xxlem7.4} Let $T$ be a finite dimensional
simple left $H$-module. Then $E(T)$ is an infinite
dimensional uniserial module and the injective
resolution of $T$ is
$$0\to T\to E(T)\to E(T_0)\to 0$$
where $T_0$ is the unique finite dimensional
simple $H$-module such that $T_0'\cong T$. If $T$ 
is 1-dimensional, then $T\cong \inth^r\otimes T_0$.
\end{lemma}

\begin{proof} Recall that $H$ is regular of global
dimension 1. Then the minimal injective resolution
of $T$ is
$$0\to T\to E(T)\to E(M)\to 0$$
for some module $M$. Since $E(T)$ is locally finite,
$E(M)$ is locally finite. By Theorem \ref{xxthm3.3}(b),
there is a unique simple module $W$ such that
$\Ext^1_H(W,T)\neq 0$ and that $W':=\Ext^1_H(W,H)^*$
is isomorphic to $T$. By notation $W=T_0$. By
the definition of $\Ext^1_H(W,T)$, $T_0$ is the 
unique simple module in the socle of $E(M)$. 
So we have the minimal injective resolution of $T$
$$0\to T\to E(T)\to E(T_0)\to 0.$$
Since $T$ and $T_0$ have no fundamental
difference, $E(T)$ is infinite dimensional over
$k$. If $N$ is a submodule of $E(T)$ of length 2,
then $\Ext^1_H(N/T,T)\neq 0$. By what we just proved,
$N/T\cong T_0$. Since $\Ext^1_H(T_0,T)=k$, $N$ is unique.

Next we use induction to show that a submodule
of $E(T)$ is uniquely determined by its length. First
we assume the length of a submodule $M\subset E(T)$,
denoted by $l(M)$, is finite. If $l(M)=1$, then $M=T$,
which is the socle of $E(T)$. Suppose the submodule of
$E(T)$ of length $<s$ is unique. Now let $M$
be a submodule of length $s$. Let $N=M/T$ be the
submodule in $E(T_0)$, which is of length $s-1$.
Since $N$ is unique in $E(T_0)$, $M$ is unique in $E(T)$.
This takes care of the case of finite dimensional
submodules $M$. Now assume that $M$ is 
infinite dimensional over $k$. Since $E(T)$
is locally finite, $M$ is not noetherian.
One can easily construct an ascending sequence
of submodules $M_s$ of length $s$. By what we
just proved, the submodules $M_s$ are uniquely determined
in $E(T)$. This means that $\{M_s\}$ is the
complete set of finite dimensional submodules.
Hence
$$E(T)=\bigcup_s M_s\subseteq M\subseteq E(T).$$

If $T$ is 1-dimensional left $H$-module,
then $T_0= \inth^r\otimes T$ since
$\Ext^1_H(\inth^r\otimes T,T) \cong \Ext^1_H(\inth^r,k)=k$.
\end{proof}

We collect and fix some notations for the rest of this
section.

\begin{convention}
\label{xxcon7.5} Let $H$ be as in Theorem \ref{xxthm7.1} 
and let $n=io(H)$.
\begin{enumerate}
\item
For any simple left $H$-module $T$, let $E(T)$ be the 
injective hull of $T$ and let $E_s(T)$ be the unique
submodule of $E(T)$ of length $s$. For every $s\geq 1$, 
let $E_s=\oplus_{i=0}^{n-1} E_s((\inth^r)^{\otimes i})$.
Let $E=\oplus_{i=0}^{n-1} E((\inth^r)^{\otimes i})$.
\item
Let $\fm_i=l.ann((\inth^r)^{\otimes i})$ and let 
$J=J_{iq}=\bigcap_{i=0}^{n-1} \fm_i$. This is the defining 
ideal of $H_{iq}$. For every $s\geq 1$, let $H_s=H/J^s$. 
Let $\widehat{H}$ be the completion $\lim_s H_s$.
Let $\gr H$ be the associated graded ring $\oplus_s J^s/J^{s+1}$.
\item
Let $\Sigma_i$ be the quotient map from $H\to H/\fm_i$. 
We identify $\Sigma_i$ with $(\Sigma^r)^i$ in $H^\circ$.
Let $\sigma_i$ be the algebra automorphism of $H$ determined
by $\Sigma_i$. Then $\sigma_i=(\sigma^r)^i$ where 
$\sigma^r=\sigma_1$. The definitions of $\Sigma^r$ and $\sigma^r$
are given before Lemma \ref{xxlem2.3}. Let $G$ be the abelian 
group $\{(\Sigma^r)^{i}\}$ (which is isomorphic to ${\mathbb Z}_n$).
Then $H_{iq}$ is isomorphic to $(kG)^\circ$.
\end{enumerate}
\end{convention}

\begin{lemma}
\label{xxlem7.6}
Assume the notation as in Convention \ref{xxcon7.5}.
\begin{enumerate}
\item
$\sigma_j(\fm_i)=\fm_{i-j}$ for all $i,j$. Equivalently,
$\sigma_j((\inth^r)^{\otimes i})=(\inth^r)^{\otimes (i-j)}$
as $\sigma_j$ acting on $H/J$.
\item
$H^{co H_{iq}}=H^{kG}$ where $G$ acts on $H$ via $\Sigma_i\to \sigma_i$.
\item
The canonical map $H\to \widehat{H}$ is injective and
every automorphism $\sigma_i$ of $H$ extends to a
automorphism of $\widehat{H}$, still denoted by $\sigma_i$.
\item
$E_s$ is the injective hull, as well as the projective cover, of
$E_1$ as left $H_s$-module.
\item
$H_s\cong E_s$ as left $H$-modules.
\end{enumerate}
\end{lemma}

\begin{proof} (a) This is equivalent to the fact 
$\pi_i
\sigma_j=\pi_i\circ \pi_j=\pi_{i+j}$ for all $i,j$.

(b) This is true because the right $H_{iq}$-coaction on
$H$ is equivalent to the right $H_{iq}^\circ$-action on $H$.

(c) The kernel of the map $H\to \lim_s H_s$ is $\bigcap_i J^i$ which
is zero by Lemma \ref{xxlem6.9}. Hence $H\to \widehat{H}$ is
injective.

By (a) $\sigma_i(J)=J$. Hence $\sigma_i(J^t)=J^t$ for all $t$.
This means that $\sigma_i$ are automorphisms of $H_s$.
Since $\widehat{H}=\lim_s H_s$, $\sigma_i$ extends to automorphisms
of $\widehat{H}$ naturally.

(d) Since $E_s(T)$ is uniserial, so is $E_s(T)/JE_s(T)$.
Thus $E_s(T)/JE_s(T)$ is simple. Thus shows that $\Hom_H(H_s,E)
=E_s$. Hence $E_s$ is the injective hull of $E_1$ as
$H_s$-module. By the analogous statement for right modules,
we see that the uniserial module $E_s^*$ is the  injective
hull of the right $H_s$-module $(E_s/{E_{s-1}})^*$.
 By $k$-linear
duality $(-)^*$, $E_s$ is the projective cover of
the $H_s$-module $E_1$, since $E_s/{E_{s-1}} \cong E_1$.

(e) The assertion follows from (d) and the fact that $H_s$ is
the projective cover of $H_s/JH_s\cong E_1$.
\end{proof}

\begin{proposition}
\label{xxprop7.7}
Assume the notation as in Convention \ref{xxcon7.5}.
\begin{enumerate}
\item
$\widehat{H}$ is isomorphic to the matrix
algebra
$$H_{M}:=\begin{pmatrix} R&xR&\cdots 
&x^{n-2}R&x^{n-1}R\\
x^{n-1}R &R &\cdots&x^{n-3}R& x^{n-2}R\\
\ldots&\ldots&\cdots& \cdots& \ldots\\
x^2R&x^3R&\cdots &R&xR\\
xR&x^2R&\cdots &x^{n-1}R&R\end{pmatrix}$$
where $R$ is a local ring $k[[x^n]]$.
\item
$\widehat{H}$ is a semilocal noetherian prime PI ring of
global dimension and Krull dimension 1.
\item
Let $e_i=H/\fm_i$ for $i=1,\cdots, n$. Then
$H/J=\widehat{H}/J(\widehat{H}) =\oplus_i ke_i$
 where $J(\widehat{H})$ is the
Jacobson radical of $\widehat{H}$.
The induced automorphism $\sigma_j$ maps $e_i$ to $e_{i-j}$.
\item
$\gr H$ is isomorphic to the algebra similar to $H_M$ where
$R$ is $k[x^n]$ instead of $k[[x^n]]$. The induced
automorphism $\sigma_j$ maps $e_i$ to $e_{i-j}$.
\item
The subring $\widehat{H}^{kG}$ is isomorphic to $k[[y]]$, whence it is
a commutative local domain of global and Krull dimension 1.
\item
$H^{kG}$ is a commutative domain.
\item
$H^{kG}$ is affine and $H$ is a finite module over $H^{kG}$
on both sides.
\end{enumerate}
\end{proposition}

\begin{proof} (a) For any map $f: M\to N$, we write $(m)f$ instead of
$f(m)$ so that we don't need to take opposite ring when we have
the canonical isomorphism $H\cong \End_H(_HH)$.

Let $e_{i,i}$ be the identity map of
$E((\inth^r)^{\oplus i})$. Let $x_{0,1}$ be the fixed
quotient morphism $E(k) \to E(\inth^r)$ in the injective
resolution of $k$. Then it induces a sequence
of quotient morphism
$$x_{i,i+1}: E((\inth^r)^{\otimes i}) \to E((\inth^r)^{\otimes (i+1)})$$
with the kernel $(\inth^r)^{\otimes i}$. Clearly,
$x_{i,i+1}=x_{i+n,i+1+n}$ for all $i$. We claim that
\begin{enumerate}
\item[(i)]
$\End_{H}(E((\inth^r)^{\otimes i}))= k[[y_i]]:=R_i$
where $y_i=x_{i,i+1}x_{i+1,i+2}\cdots x_{i+n-1,i+n}$ for all
$i\in {\mathbb Z}_n$.
\item[(ii)]
For $0<j-i<n$, $\Hom_H(E_s((\inth^r)^{\otimes i}),
E_s((\inth^r)^{\otimes j}))=x_{i,i+1}\cdots x_{j-1,j}R_j=
R_i x_{i,i+1}\cdots x_{j-1,j}$.
\end{enumerate}
To prove (i) we note that $\End_{H}(E_n((\inth^r)^{\otimes i}))=k$.
By induction on $s$ one shows that
$\End_{H}(E_{sn}((\inth^r)^{\otimes i}))=k[y_i]/(y_i^s)$.
Now we have
$$\End_{H}(E((\inth^r)^{\otimes i}))=
\lim_s \End_{H}(E_{sn}((\inth^r)^{\otimes i}))=k[[y_i]].$$
The proof of (ii) is similar.

Identifying $x_{i,i+1}$ with $e_{i,i+1}x$ inside the matrix
algebra $M_{n}(k[[x]])$ where $\{e_{i,j}\}$ is the matrix unit,
we have an isomorphism $\End_H(E)\cong H_M$. By Lemma
\ref{xxlem7.6}(e), $H_s\cong E_s$ as left $H$-modules.
So we can identify $E$ with $\lim_s (H_s)^*$. Hence
$$\widehat{H}=\lim_s H_s=\lim_s \End_H(H_s)\cong
\lim \End_H(E_s)=\End_H(E).$$
Therefore (a) is proved.

(b) This is clear.

(c) This follows from Lemma \ref{xxlem7.6}(a) and the fact
$H/J=\widehat{H}/J(\widehat{H})$.

(d) Since $\gr H=\gr \widehat{H}\cong \gr H_M$.
This is clear from the description of (a). Equivalently,
in the notations in the proof of (a), $\gr H$
is isomorphic to
$$k\langle e_{0,0},\cdots, e_{n-1,n-1}, x_{0,1},x_{1,2},\cdots,x_{n-1,0}
\rangle/(rels)$$ where the ideal $(rels)$ is generated by the relations
$1=\sum_{i=0}^{n-1}e_{i,i}$,
$e_{i,i} x_{i,k}=x_{i,k}=x_{i,k}e_{k,k}$, and $0=e_{i,i} x_{j,k}=
x_{k,j}e_{i,i}=x_{i,j}x_{k,l}$ if $j\neq k$ in ${\mathbb Z}_n$.

Since $\sigma_j$ maps $J$ to $J$, $\sigma_j$ induces
an automorphism of $\gr H$. Therefore
$kG\cong \oplus_i k\sigma_j$ acts on $\gr H$.
The second assertion follows from (c).

(e) First we compute the the invariant subring $(\gr H)^{kG}$
of $\gr H$. By (d), $\sigma:=\sigma_{-1}$ maps $e_{i,i}$ to
$e_{i+1,i+1}$. Use the equation
$$\sigma(x_{i,i+1})=\sigma(e_{i,i}x_{i,i+1}e_{i+1,i+1})
=e_{i+1,i+1}\sigma(x_{i,i+1})e_{i+2,i+2}$$
one sees that $\sigma(x_{i,i+1})=w_i x_{i+1,i+2}$ for some
$0\neq w_i\in k$. Since $\sigma^n$ is the identity,
we obtain that $w_0\cdots w_{n-1}=1$. Modifying $x_{i,i+1}$ by
scalars, we may assume that $w_i=1$ for all $i$. Now $\sigma$
maps $e_{i,i}\to e_{i+1,i+1}$ and $x_{i,i+1}\to x_{i+1,i+2}$
for all $i\in {\mathbb Z}_n$. Now it is straightforward
to compute that $(\gr H)^{kG}=k[x]$ where
$x=\sum_{i=0}^{n-1}x_{i,i+1}$.

We now study the subring $(\widehat{H})^{kG}$.
Let $y=\sum_{i=0}^{n-1} \sigma^i(x_{i,i+1})$. Then $\gr y=x$
because $\gr y$ is a $\sigma$ invariant and inside
$\gr H$ one has $\sum_{i=0}^{n-1} \sigma^i(x_{i,i+1})=x$.
We will show that $(\widehat{H})^{kG}=k[[y]]$. Define
the degree of $f\in \widehat{H}$ to be the degree $\gr f$
in the $\gr H$. We claim that for every element $f$ in
$(\widehat{H})^{kG}$, there is a $w\in k$ such that
$\deg(f-w y^p)>p$ where $p=\deg f$. This is true because
$\gr f$ is in $(\gr H)^{kG}$ and hence $\gr f$ and $\gr y^p$
are in the same vector space $kx^p$. So we proved the
claim. By induction and the claim, for any element $f\in
\widehat{H}^{kG}$, there is a sequence of elements $w_i\in k$ such that
$\deg (f-\sum_{i=0}^{t} w_i y^i)> t$ for all $t$. Since
$\widehat{H}$ is complete $f=\sum_{i=0}^{\infty} w_i y^i
\in k[[y]]$.

(f) Since $H\subset \widehat{H}$ we have $H^{kG}\subset
(\widehat{H})^{kG}$. By (e) $H^{kG}$ is a commutative domain.

(g) [Brown] Let $Z$ be the center of $H$. Since $H$ is affine 
prime of GK-dimension one, $Z$ is an affine commutative 
domain and $H$ is a finite module over $Z$ on both sides
\cite{SSW}. 
Since $G$ is a finite group of automorphism of $H$, it is a 
group of automorphism of $Z$. Since $Z$ is affine, 
the invariant subring $Z^G$ is affine and $Z$ is a finite
module over $Z^G$ by Hilbert-Noether theorem 
\cite[Theorem 1.3.1]{Be}. Hence $H$ is a finite module
over $Z^G$. Since $H^{kG}$ is a subring of $H$ containing 
$Z^G$, $H^{kG}$ is a finite module over an affine 
commutative subring $Z^G$. Therefore $H^{kG}$ is affine 
and $H$ is a finite module over $H^{kG}$ on both sides. 

As a consequence, $H^{kG}$ has GK-dimension and Krull
dimension 1.
\end{proof}

\begin{proof}[Proof of Theorem \ref{xxthm7.1}]
(a) If $H$ is unimodular, then, by Proposition
\ref{xxprop7.7}, $H$ is a subring of $\widehat{H}=
M_1(k[[x]])$. Hence $H$ is commutative. 

For the second assertion we assume that  $k$ is 
algebraically closed. The algebraic group 
associated to $H$ is connected, linear (i.e., 
affine) and of dimension 1. By 
\cite[Theorem 2.6.6]{Sp}, any connected 
linear algebraic group of dimension 1 is isomorphic
to either ${\mathbb G}_a:=(k,+)$ or ${\mathbb G}_{m}:=
(k-\{0\},\times)$. Thus $H$ is isomorphic to either 
$k[x]$ or $k[x^{\pm 1}]$.

(b) Suppose $n=io(H)>1$. By Proposition \ref{xxprop7.7}(a),
the PI degree of $\widehat{H}$ is no more than $n$.
Hence the PI degree of $H$ is no more than $n$.
Combining with Lemma \ref{xxlem5.3}(g), $io(H)=PI.deg(H)$.
This proves (i). The rest of (ii,iii) follows from Proposition
\ref{xxprop7.7}(f,g). 
\end{proof}

\begin{corollary}
\label{xxcor7.8}
Let $H$ be as in Theorem \ref{xxthm7.1}. 
Suppose $k$ is algebraically closed.
The following are equivalent:
\begin{enumerate}
\item
$H$ is a domain.
\item
$H$ is commutative.
\item
$\epsilon(\inth^l)\neq 0$ where $\epsilon$ is defined in 
\eqref{E3.0.1} (for $d=1$).
\item
$H$ is unimodular.
\item
$\Ext^1_H(k,k)\neq 0$.
\item
$H$ is isomorphic to either $k[x]$ or $k[x^{\pm 1}]$.
\end{enumerate}
\end{corollary}

\begin{proof} (a) $\Rightarrow$ (b) 
Let $Q$ be the ring of fractions of 
$H$. Hence $Q$ is a division algebra 
of GK-dimension 1. We claim that every 
division algebra of GK-dimension 1 
is commutative, which implies (b).
It is clear that $Q$ is a direct union
of finitely generated division algebras.
We may assume that both $Q$ and 
its center, say $K$, are finitely 
generated as a division algebras. 
Since the transcendence degree of $K$
is 1 and the base field $k$ is algebraically
closed, by Tsen's theorem \cite[P.~374]{Co}, 
the Brauer group of $K$ is trivial, whence 
$Q=K$. Therefore $Q$ is commutative. 

(b) $\Rightarrow$ (c) Since $H$ is commutative, 
$\inth^r=\inth^l=k$, (c) follows from (Cond1)
since $H$ is regular. 

(c) $\Rightarrow$ (d) This is Lemma \ref{xxlem3.1}. 

(d) $\Rightarrow$ (e) 
$\Ext^1_H(k,k) \cong \Ext^1_H(\inth^r,k) \cong k$.

(e) $\Rightarrow$ (c) By (Cond3) for $T=k$, we see
that $\epsilon(\inth^l)\neq 0$. 

(d) $\Rightarrow$ (f) This is Theorem \ref{xxthm7.1}(a).

(f) $\Rightarrow$ (a) Trivial.
\end{proof}

\begin{lemma}
\label{xxlem7.9}
Let $H$ be a noetherian regular affine Hopf
algebra of GK-dimension 1. Let $W$ and $T$ 
be left simple $H$-modules in the same clique.
Then $\dim W=\dim T$, $\End_H(W)\cong \End_H(T)$
and $PI.deg(H/l.ann_H(W))=PI.deg(H/l.ann_H(T))$.
\end{lemma}

\begin{proof} First of all the global
dimension of $H$ will be 1 [Lemma \ref{xxlem5.3}(a)].
By Theorem \ref{xxthm3.3}, for every given 
simple module $S$, there is 
only one simple module $V=S'$ such that 
$\Ext^1_H(S,V)\neq 0$. In this case we write
$V\leadsto S$. This uniqueness property implies 
there is a unique sequence of simple modules
$\{T_0, T_1,\cdots T_n\}$ such that 
either
$$T=T_0\leadsto T_1\leadsto \cdots \leadsto T_n=W$$
or
$$W=T_0\leadsto T_1\leadsto \cdots \leadsto T_n=T.$$
By induction on $n$ and the left-right symmetry we 
may assume that 
$T\leadsto W$. In this case $T=W'=\Ext^1_H(W,H)^*$.
Since $H$ is AS-Gorenstein, $\Ext^1_H(-,H)$ is 
an equivalent functor $\Mod_{fd}$-$H\to \Mod_{fd}$-$H^{op}$.
Hence we have $\End_H(W')=\End_H(W)$. By (AS2)$'$ in
Section 1, $\dim W'=\dim W$. Finally $PI.deg(H/l.ann_H(W)$ 
is determined by $\dim W$ and $\dim \End_H(W)$.
\end{proof}

\begin{proof}[Proof of Corollary \ref{xxcor0.3}]
By \cite[Remark 5.7(ii)]{SZ}, two
maximal ideals $I$ and $J$ of $H$ are in the 
same clique if and only if $I\cap Z(H)=J\cap Z(H)$
where $Z(H)$ is the center of $H$. Let $I=l.ann_H(M)$.
By \cite[Theorem 8]{Bra} we have 
\begin{equation}
\label{E7.10.1}
PI.deg(H)=\sum_{J} c_J PI.deg(H/J),
\tag{E7.10.1}
\end{equation}
for some integers $c_J$, where $J$ runs over all 
maximal ideals in the clique containing $I$.
Let $M_J$ be the simple left 
$H$-modules corresponding to $J$. By Lemma 
\ref{xxlem7.9} $\dim M_J=\dim M$ and
$PI.deg(H/I)=PI.deg(H/J)$ for all $J$. 
Since $k$ is algebraically closed, 
$PI.deg(H/I)=\dim M$. Then \eqref{E7.10.1}
becomes 
$$PI.deg(H)=(\sum_{J} c_J) \dim M.$$
The assertion follows by Theorem 
\ref{xxthm0.2}(a).
\end{proof}

Note that we can show that coefficients $c_J$'s
in \eqref{E7.10.1} are 1, but the proof is omitted.

\section{Group algebras}
\label{xxsec8}

There is nothing new in this section. The reason we 
include this short section is to show how 
integrals and integral order can be related to the 
structures of groups. The following lemma is well-known.

\begin{lemma}
\label{xxlem8.1} Let $G$ be a group and $H'$ be a 
Hopf algebra. Suppose $f: kG\to H'$ is a surjective
Hopf algebra homomorphism. Then 
\begin{enumerate}
\item 
$H' = kG'$, where $G'=f(G)$, and $f$ is induced 
by the  group homomorphism $f|_G: G\to G'$. 
\item
$kG_0=(kG)^{co\; kG'}$, where $G_0$ is the 
kernel of $f|_G$.
\end{enumerate}
\end{lemma}

Let $G$ be a finitely generated group. We say 
$G$ has linear growth if (a) $G$ is infinite and
(b) there is a generating set $T\subset G$ 
with $1\in T$ and $T^{-1}=T$ such that $|T^n|\leq 
cn$ for some constant $c>0$. By \cite{WV}, 
if $G$ is sub-quadratic (meaning that 
$|T^n|-|T^{n-1}|<n$ for some $n$, then 
$G$ has linear growth and contains a 
subgroup ${\mathbb Z}$ of finite index. 
As a consequence, $kG$ is noetherian of GK-dimension
1. Such a group $G$ is said to have 
{\it two ends} \cite[p.~100]{IS}. 
Conversely, if $kG$ is a noetherian Hopf
algebra of GK-dimension 1, then $G$ is finitely
generated with linear growth.
By \cite[Theorem 3.13 in Chapter 10]{Pas} $kG$ is regular if 
$k$ is of characteristic zero (this can be weakened, 
dependent on the order of finite subgroups of $G$).

The conjectural descriptions \eqref{E0.2.1} and 
\eqref{E0.2.2} are verified for group rings. 
The following result is basically known \cite{IS,St}
and we remark that it gives the short exact sequences
\eqref{E0.2.1} and \eqref{E0.2.2}.

\begin{proposition}
\label{xxprop8.2} 
Let $G$ be a finitely generated group with linear 
growth. 
\begin{enumerate} 
\item 
There is a finite normal subgroup $G_{dis}\subset G$
such that it contains all finite normal subgroups 
of $G$. And there is a short exact sequence
$$1\to G_{dis}\to G\to G_{conn}\to 1$$
where $G_{conn}=G/G_{dis}$ contains no nontrivial
finite normal subgroup.
\item
The connected component $G_{conn}$ is isomorphic 
to either ${\mathbb Z}$ or ${\mathbb D}$ [Example
\ref{xxex4.6}]. In either
case, $G_{conn}$ fits into a short exact sequence 
$$1\to {\mathbb Z}\to G_{conn}\to G_{iq}\to 1$$
where $G_{iq}$ is either $\{1\}$ or ${\mathbb Z}_2$.
\item
$io(kG)=1$ or $2$, corresponding to the two cases in (b).
\end{enumerate}
\end{proposition}

It follows from Proposition \ref{xxprop8.2}(a,b) that 
\eqref{E0.2.1} and \eqref{E0.2.2} hold for $H=kG$ when
$kG$ is affine regular of GK-dimension 1.

\begin{proof}[Proof of Proposition \ref{xxprop8.2}] 
(a) We take $k={\mathbb C}$ for simplicity
and let $H=kG$. Then $kG$ is a regular noetherian affine 
PI Hopf algebra of GK-dimension 1. By Theorem 
\ref{xxthm6.5}, there is a surjective Hopf algebra 
homomorphism $f: kG\to H_{coon}$ with $H_{conn}$ being
prime. By Lemma \ref{xxlem8.1}(a), $H_{conn}$ is a 
group algebra, denoted by $kG_{conn}$, and 
$G_{conn}$ is a quotient group of $G$. Since $kG_{conn}$
is prime, $G_{conn}$ contains no non-trivial finite 
normal subgroup. We obtain the short exact sequence 
by letting $G_{dis}=\ker(G\to G_{conn})$. Since 
$G_{conn}$ is infinite and $G$ has linear growth, 
$G_{dis}$ is finite. By Theorem \ref{xxthm6.5}(b), 
$G_{dis}$ contains all finite normal subgroups of $G$. 

(b) Now assume that $kG$ is prime. If $kG$ is commutative,
then $G$ is abelian and it must be ${\mathbb Z}$ by the
decomposition of finitely generated abelian groups.
Now assume $kG$ is not commutative. By Theorem 
\ref{xxthm7.1}, there is a Hopf algebra homomorphism
$kG\to H_{iq}$ such that $(kG)^{co \; H_{iq}}$ is a 
commutative domain. Note that $H_{iq}=
(k{\mathbb Z}_n)^\circ\cong k{\mathbb Z}_n$ 
where $n=io(kG)$. By Lemma \ref{xxlem8.1}, 
$(kG)^{co \; H_{iq}}$ is a group algebra of 
GK-dimension 1, hence $(kG)^{co \; H_{iq}}=k{\mathbb Z}$.
So we have a short exact sequence 
\begin{equation}
\label{E8.2.1}
1\to {\mathbb Z}\to G\to {\mathbb Z}_n\to 1,
\tag{E8.2.1}
\end{equation}
which gives rise to the canonical map $kG\to 
(kG)_{iq}$. Let $K$ be a maximal normal abelian 
subgroup of $G$ containing ${\mathbb Z}$ in 
\eqref{E8.2.1}. So we have a short exact sequence
$$1\to K\to G\to {\mathbb Z}_m\to 1$$
for $1< m\leq n$. Since $G$ does not contain
proper finite normal subgroup, it is easy to 
show that $K$ does not contain proper finite 
subgroup. Thus $K\cong {\mathbb Z}$. If the 
action of ${\mathbb Z}_m$ on $K$ is not
faithful, then it produces a larger normal 
abelian group containing $K$, a contradiction.
Therefore the action of ${\mathbb Z}_m$ on 
$K$ is faithful. The only non-trivial action 
of $K\cong {\mathbb Z}$ is $n\to -n$. 
Thus $m=2$. In this case $G\cong {\mathbb D}$,
which is described in Example \ref{xxex4.6}.

(c) Follows from (b) and Example \ref{xxex4.6}.
Note that the assertion holds for any field $k$
such that $kG$ is regular. 
\end{proof}

The following example is provided by Stafford.

\begin{example}
\label{xxex8.3} Let $k$ be a field of 
$\operatorname{char}k\neq 2$. Suppose that 
$i:=\sqrt{-1}\not\in k$. One might assume
$k={\mathbb R}$ for simplicity. Let $H$ be the
algebra $k[x,y]/(x^2+y^2-1)$. This is a Hopf algebra
because it is the coordinate ring of the unit
circle, which is an algebraic group at least when
$k={\mathbb R}$. The coalgebra 
structure of $H$ is determined by 
$$\Delta(x)=x\otimes x-y\otimes y, 
\quad \Delta(y)=x\otimes y+y\otimes x$$
and the counit and the antipode are 
determined by
$$\epsilon(x)=1, \epsilon(y)=0\quad \text{and}\quad
S(x)=x, S(y)=-y.$$
Since $i\not\in k$, this Hopf algebra is not
isomorphic to the group algebra $k {\mathbb Z}
(\cong k[t,t^{-1}])$ over the base field $k$.
Let $F$ be any field extension of $k$ such that
$i\in F$, or one can take $F=\bar{k}$.
Then $H\otimes_k F$ is isomorphic to $F[z,z^{-1}]
\cong F{\mathbb Z}$ where $z=x+iy$ and $z^{-1}=x-iy$. 
So there are two non-isomorphic Hopf algebras, namely,
$k{\mathbb Z}$ and $H$, such that their
field extensions are isomorphic as Hopf algebras. 

A Hopf quotient of $H$ was studied by
Greither and Pareigis \cite{GP}. Let
$H_1=H/(xy)$. Then $H_1$ is a finite dimensional
Hopf algebra not isomorphic to a group
algebra, and $H_1\otimes_k F$ is isomorphic to 
a group algebra 
$F ({\mathbb Z}_2\times {\mathbb Z}_2)$. This
Hopf algebra is called the {\it circle
Hopf algebra} or the {\it trigonometric Hopf 
algebra} (see also \cite{Par} and \cite[pp.125-6]{Mo}).

An easy extension can be made without 
restriction on $\operatorname{char} k$ as follows. Let 
$\xi\in k$ such that $\sqrt{-\xi}\not\in k$. Let 
$H_{\xi}$ be the algebra $k[x,y]/(x^2+\xi y^2-1)$ 
with other operations
$$\Delta(x)=x\otimes x-\xi y\otimes y, 
\quad \Delta(y)=x\otimes y+y\otimes x$$
and 
$$\epsilon(x)=1, \epsilon(y)=0\quad \text{and}\quad
S(x)=x, S(y)=-y.$$
Then $H_{\xi}$ is a Hopf algebra.

It is easy to check that if $K$ is a Hopf algebra
such that $K\otimes_k F$ is isomorphic to the Hopf
algebra $F[z]$ for some field extension $F\supset k$, 
then $H\cong k[x]$.

One can also construct a similar Hopf algebra $H'$
such that $H'\not\cong k {\mathbb D}$ 
(see Example \ref{xxex4.6} for the definition of
${\mathbb D}$); but for some field extension 
$F$, $H'\otimes_k F\cong F {\mathbb D}\cong
k {\mathbb D}\otimes_k F$.
\end{example}

\begin{question}
\label{xxque8.4} When $k$ is algebraically closed,
the group algebras in Proposition 
\ref{xxprop8.2} and the Hopf
algebras in Example \ref{xxex2.7}
(for $m=1$) are the only examples of 
affine prime regular Hopf algebras of 
GK-dimension 1 we know so far. Are there
others?

When $k$ is not algebraically closed, 
there are some others [Example \ref{xxex8.3}].
What can we expect in this case?
\end{question}

Finally we give an example that shows that Theorem 
\ref{xxthm0.2} fails for group algebras of 
GK-dimension 2.

\begin{example}
\label{xxex8.5}
Let $G$ be the group $\langle g,x,y |
g^2=1, gxg=x^{-1},gyg=y^{-1}, xy=yx\rangle$. 
It's easy to see that the subgroup $\langle y\rangle$ 
is normal and $G/\langle y\rangle $ is isomorphic to
the group ${\mathbb D}$ defined in Example \ref{xxex4.6}.

By using the relations between $g,x$ and $y$, we see that
$$g(y-1)=(y^{-1}-1)g=(y-1)(-y^{-1}g), \quad
\text{and}\quad x(y-1)=(y-1)x.$$
Hence $y-1$ is a normal element in the group algebra
$kG$ and $kG/(y-1)\cong k{\mathbb D}$. By Example 
\ref{xxex4.6} the right integral of $k{\mathbb D}$ 
is isomorphic to $k{\mathbb D}/(x-1,g+1)$. By Lemma 
\ref{xxlem2.6}(b), the right integral of $kG$ is 
isomorphic to the trivial $kG$-module $k$. Hence 
the integral order is 1. This shows the first 
property listed below. The other properties are clear. 
\begin{enumerate}
\item
$io(kG)=1$.
\item
$kG$ is a prime ring.
\item
$kG$ has global dimension and GK-dimension 2.
\item
$kG$ is not a domain.
\item
$kG$ is a PI ring of PI degree 2, hence not commutative.
\end{enumerate}
Therefore all statements in Theorem \ref{xxthm0.2} fail 
for this group algebra.
\end{example}

\section*{Acknowledgments}

Q. -S. Wu is supported by the NSFC (key project 10331030) 
and  by  STCSM (03JC14013) and supported by the 
Cultivation Fund of the Key Scientific and Technical 
Innovation Project, Ministry of Education of China 
(NO 704004). J.J. Zhang is supported 
by NSF grant DMS-0245420 (USA) and Leverhulme 
Research Interchange Grant F/00158/X (UK). The 
authors thank Jacques Alev, Ken Brown, Ken Goodearl, 
Tom Lenagan, Martin Lorenz, Monty McGovern, Susan Montgomery, 
John Palmieri, Paul Smith and Toby Stafford for 
many useful discussions and valuable comments.
In particular the authors thank Toby Stafford for 
providing Example \ref{xxex8.3} and thank Ken Brown 
for the proofs of Lemma \ref{xxlem6.8}(a) and 
Theorem \ref{xxthm7.1}(b)(iii).

\providecommand{\bysame}{\leavevmode\hbox
to3em{\hrulefill}\thinspace}
\providecommand{\MR}{\relax\ifhmode\unskip\space\fi MR
}
\providecommand{\MRhref}[2]{%

\href{http://www.ams.org/mathscinet-getitem?mr=#1}{#2}
}
\providecommand{\href}[2]{#2}

\end{document}